\documentclass{svjour3}

\RequirePackage{fix-cm}
\usepackage{graphicx}
\usepackage{psfrag}
\usepackage{subfigure}
\usepackage{url}
\usepackage{color}
\usepackage{cite}
\usepackage{epsfig}
\usepackage{hyperref}
\usepackage{amsmath,amssymb}
\usepackage{slashbox}
\usepackage{tablefootnote}
\usepackage{longtable}

\newcommand{\rset}{\mathbb{R}}

\newcommand{\Ncal}{\mathcal{N}}
\newcommand{\Scal}{\mathcal{S}}

\newcommand{\lu}{\mathcal{L}}
\newtheorem{assumption}[theorem]{Assumption}

\usecounter{algorithmenumi}

\date{October 2013}

\begin{document}

\title{Distributed dual gradient methods and error bound conditions}

\author{Ion Necoara and Valentin Nedelcu
\thanks{The authors are with Automatic Control and Systems Engineering
Department, University Politehnica Bucharest, Romania.}}

\titlerunning{Converge rate of distributed dual first order methods}

\maketitle

\begin{abstract}
In this paper\footnote{This paper is based on Chapter 4 of the Ph. D. Thesis \cite{Ned:13}.} we propose distributed dual gradient algorithms for
linearly constrained separable convex problems and analyze their
rate of convergence under different assumptions. Under the strong
convexity assumption on the primal objective function we propose
two distributed dual fast gradient schemes for which we prove
sublinear rate of convergence for dual suboptimality but also
primal suboptimality and feasibility violation for an average
primal sequence or for the last generated primal iterate.  Under
the additional assumption of Lipshitz continuity of the gradient
of the primal objective function we prove a global error bound
type property for the dual problem and then  we analyze a dual
gradient scheme for which we derive  global linear rate of
convergence for both dual and primal suboptimality and primal
feasibility violation. We also provide numerical
simulations on optimal power flow problems.
\end{abstract}

\section{Introduction}
\label{sec_introduction}
Nowadays, many engineering applications which appear in the context
of communications networks or networked systems can be posed as
linearly constrained  separable convex problems.  Several important
applications that can be modeled in this framework,  the network
utility maximization (NUM) problem  \cite{BecNed:13}, the optimal power
flow (DC-OPF) problem for a power system \cite{ZimMur:11} and distributed
model predictive control (MPC) problem  for networked systems
\cite{NecNed:13},  have attracted great attention lately. Due to the
large dimension and the separable structure of these problems,
distributed optimization methods have become an appropriate tool for
solving such problems.

\noindent The standard approach to distributed optimization in networks is to
use  decomposition.  Decomposition methods represent a powerful tool
for solving these type of problems due to their ability of dividing
the original large scale problem into smaller subproblems which are
coordinated by a master problem. Decomposition methods can be
divided in two main classes: primal and dual decomposition. While in
the primal decomposition methods the optimization problem is solved
using the original formulation and variables, in dual decomposition
the constraints are moved into the cost using the Lagrange
multipliers and the dual problem is solved. In many applications,
such as (NUM), (DC-OPF) and (MPC) problems, when the constraints set is
complicated (i.e. the projection on this set is hard to compute)
dual decomposition becomes more effective since a primal approach
will require at each iteration a projection onto the feasible set,
operation that is numerically very expensive.

\noindent First order methods for solving dual problems  have been
extensively studied in the literature. Subgradient methods based
on averaging, that produce primal solutions in the limit, can be
found e.g. in \cite{KiwLar:07, LarPat:98, SenShe:86}. Despite
widespread use of the (sub)gradient methods for solving dual
problems, there are some aspects that have not been fully studied.
In particular, in practical applications, the main interest is in
finding an approximate primal solution that is near-feasible and
near-optimal. Moreover, we need to characterize the convergence
rate for the approximate primal solution.  Finally, we are
interested in providing distributed schemes, i.e. methods based on
distributed computations. These represent the main issues that we
pursue in this paper.

\noindent Convergence rate analysis for the dual subgradient method  has been
studied e.g. in \cite{NedOzd:09,DoaKev:11}, where estimates of order
$\mathcal{O} (1/\sqrt{k})$ for suboptimality and  feasibility
violation of an average primal sequence are provided, with $k$
denoting the iteration counter. In \cite{NecSuy:08} the authors
propose a dual fast gradient algorithm based on a smoothing
technique and prove  rate of convergence of order
$\mathcal{O}\left(\frac{1}{k}\right)$ for primal suboptimality and
feasibility violation for an average primal sequence.     Also, in
\cite{NecNed:13} the authors propose  inexact dual (fast) gradient
algorithms for which estimates of order
$\mathcal{O}\left(\frac{1}{k}\right)$
 ($\mathcal{O}\left(\frac{1}{k^2}\right)$) in an average primal
 sequence  are provided for both primal and dual suboptimality and
 primal feasibility violation. For the
special case of QPs problems, dual gradient algorithm were also
analyzed in \cite{KogFin:11,PatBem:12,RicMor:11}.  From our
knowledge first result on the  linear convergence of dual gradient
method was provided in \cite{LuoTse:93a}. However, the authors in
\cite{LuoTse:93a} were able to show linear convergence only
locally. Finally, very few results  were known in the literature
on distributed implementations of dual gradient type methods since
most of the papers enumerated above require a centralized  step
size. Recently, the authors in \cite{BecNed:13} propose a
distributed dual fast gradient algorithm where the step size is
chosen distributively and provide estimates of order
$\mathcal{O}\left(\frac{1}{k}\right)$ for primal suboptimality and
feasibility violation in the last primal iterate. All of these
 limitation motivates our work here.

\noindent In this paper we propose distributed versions of dual
first order methods generating approximate primal feasible and
primal optimal solutions but with great improvement on the
convergence rate w.r.t. the existing results from the literature. In
particular, under the strong convexity assumption on the primal
objective function we derive a distributed version of the dual fast
gradient algorithm presented in \cite{NecNed:13} for which we are
able to provide estimates of order
$\mathcal{O}\left(\frac{1}{k^2}\right)$ on primal suboptimality and
feasibility violation for an average primal sequence. In comparison
with the algorithm proposed in \cite{NecNed:13} we do not require a
centralized  step size and thus we derive a distributed
implementation of the algorithm. Also, the estimates on primal
suboptimality and feasibility violation for our distributed
algorithm are with an order of magnitude better that the ones of
algorithm given  in \cite{BecNed:13}. We also propose a hybrid dual
fast gradient algorithm which allows us to provide estimates of
order $\mathcal{O}\left(\frac{1}{k\sqrt{k}}\right)$ on primal
suboptimality and feasibility violation in the last primal iterate.
Note that also in this case the iteration complexity of our method
is better than  of the method given in \cite{BecNed:13}. Under the
additionally Lipschitz continuity assumption on the gradient of
primal objective function, which is often satisfied in practical
applications (e.g. (NUM) and (MPC) problems), we prove that the
corresponding dual problem satisfies  a certain error bound property
\cite{LuoTse:93a}. In order to prove such a property we extend the
approach developed in \cite{WanLin:13,LuoTse:93a} to the case when
the constraints set is an unbounded polyhedron. In these settings we
analyze the convergence behavior of a distributed dual gradient
algorithm for which we are able to provide for the first time
\textit{global} linear convergence rate on primal suboptimality and
infeasibility for the last primal iterate, as opposed  to the
results in \cite{LuoTse:93a} where only \textit{local} linear
convergence was derived for such an algorithm.   We also show that
the theoretical estimates on the convergence rate depend on a
natural and easily computable measure of separability of the
problem.

\noindent \textit{ Contribution.}  In summary, the contributions of
this paper include:
\begin{itemize}

\item[(i)] We propose and analyze novel  dual gradient type
algorithms having distributed implementations   and fast rate of
convergence that generate approximate primal solutions for separable
(smooth) convex problems with linear constraints.

\item[(ii)] For these distributed algorithms we derive estimates on primal
suboptimality and infeasibility  in an average/last sequence: a dual
fast gradient method with convergence rate $\mathcal{O} ( 1 / k^2)$
in an average primal sequence; an hybrid dual fast gradient method
with convergence rate $\mathcal{O}( 1 / k^{3/2})$ in the last primal
iterate; a dual  gradient method with linear convergence in the last
primal~iterate.

\item[(iii)] Under strong convexity and  Lipschitz continuity
 of the gradient of the primal objective function  we prove an error
bound type property for the dual problem which allows us to obtain
global linear convergence for a distributed dual gradient method.

\end{itemize}

\noindent \textit{Paper Outline.} In Section \ref{sec_formulation}
we introduce our optimization model and discuss several practical
applications which can be posed in this framework. In Sections
\ref{sec_ddfg} and \ref{sec_hddfg} we propose two distributed dual
fast  gradient algorithms and provide sublinear estimates for both
dual and primal suboptimality, but also for primal feasibility
violation in an average primal sequence or in the last generated primal iterate.
In Section \ref{sec_dg_error_bound} we show that under
additional assumptions on the primal objective function the dual
problem has some error bound property which allows us to prove
 global linear converge for a distributed dual gradient method. Finally, in
Section \ref{sec_numerical} we provide extensive numerical
simulations in order to certify our proposed theory.

\noindent \textit{Notations:} We work in the space $\rset^n$
composed of column vectors. For $z, y \in \rset^n$ we denote the
standard Euclidean inner product $\langle z,y \rangle= \sum_{i=1}^n
z_i y_i$, the Euclidean norm $\left \| z \right \|=\sqrt{\langle z,
z \rangle}$ and the infinity norm $\|z\|_{\infty}=\sup_i |z_i|$.
Also, w.r.t. to the Euclidean norm $\|\cdot\|$ we denote the
projection onto the non-negative orthant $\rset^n_+$ by $\left[ z
\right]_+$ and the projection onto the convex set $D$ by  $[z]_D$.
For a positive definite matrix $W$ we denote the weighted norm of a
vector $z$ by $\|z\|_W^2=z^TWz$ and the projection of the vector $z$
onto a convex set $D$ w.r.t. to norm  $\|\cdot\|_W$ by $[z]_D^W$.
For a (block) matrix $A$ we define by $A_{i}$ its $i$th (block)
column. We denote by $I_q$ the identity matrix in $\rset^{q \times
q}$ and by $0_{p,q}$ the matrix from $\rset^{p \times q}$ with all
entries zero.


\section{Problem formulation}
\label{sec_formulation}

We consider the following linearly constrained separable convex
optimization problem:
\begin{align}
\label{eq_prob_princc} f^* = &\min_{z_i \in \rset^{n_i}}
f(z)~~~\left(= \sum_{i=1}^M f_i(z_i)\right)\\
&\text{s.t.:} ~~Az=b, ~~Cz\leq c, \nonumber
\end{align}
where $f_i$ are convex functions, $z=\left[z_1^T \cdots
z_M^T\right]^T$, $A \in \rset^{p \times n}$, $C \in \rset^{q
\times n}$, $b\in \rset^p$ and $c \in \rset^q$. To our
optimization problem \eqref{eq_prob_princc} we associate a
communication bipartite graph $\mathcal{G}=\left(V_1,V_2,E\right)$, where
$V_1=\left\{1,\dots,M\right\}$,
$V_2=\left\{1,\dots,\bar{M}\right\}$ and $E \in
\left\{0,1\right\}^{\left(\bar{M}\right)\times {M}}$ is an
incidence matrix.  We also introduce the index sets
$\Ncal_i=\left\{j \in V_2 : E_{ij} \neq 0\right\}$ for all $i \in
V_1$ and $\bar{\Ncal}_j=\left\{i \in V_1 : E_{ij} \neq 0\right\}$
for all $j \in V_2$ which describe the local information flow in
the graph. Note that the cardinality of the sets $\Ncal_i$ and
$\bar{\Ncal}_j$ can be viewed as a measure for the degree of
separability of problem \eqref{eq_prob_princc}. Therefore, the
local information structure imposed by the graph $\mathcal{G}$ should be
considered as part of the problem formulation. We assume that $A$
and $C$ are block matrices with the blocks $A_{ij} \in \rset^{p_i
\times n_i}$ and $C_{ij} \in \rset^{q_i \times n_i}$, where
$\sum_{i=1}^M n_i=n$, $\sum_{i=1}^{\bar{M}} p_i=p$ and
$\sum_{i=1}^{\bar{M}} q_i=q$. We also assume that if $E_{ij}=0$,
then both blocks $A_{ij}$ and $C_{ij}$ are zero. In these settings
we allow a block $A_{ij}$ or $C_{ij}$ to be zero even if
$E_{ij}=1$.

\noindent Further, we make the following assumption on the
optimization problem \eqref{eq_prob_princc}:
\begin{assumption}
\label{ass_strong}
\begin{itemize}
\item[(a)] The functions $f_i$ are $\sigma_i$-strongly convex w.r.t. Euclidean
norm $\|\cdot\|$ \cite{Nes:04}.
\item[(b)] The feasible set of problem \eqref{eq_prob_princc} is
nonempty and there exists $\bar z$ such that $A \bar z = b$ and $C
\bar z < c$.
\end{itemize}
\end{assumption}

Note that if Assumption \ref{ass_strong} $(a)$ does not hold, we
can apply smoothing techniques by adding a regularization term to
the function $f_i$ in order to obtain a strongly convex
approximation of it (see e.g. \cite{NecSuy:08} for more details).
Assumption \ref{ass_strong} $(b)$ implies that strong duality
holds for optimization problem \eqref{eq_prob_princc} and the set
of optimal Lagrange multipliers is bounded \cite{Gau:77,Man:85}.
In particular, we have:
\begin{equation}
\label{eq_dual_prob} f^*=\max_{\nu \in \rset^p,\mu \in \rset^q_+}
d(\nu,\mu),
\end{equation}
where $d(\nu,\mu)$ denote the dual function of
\eqref{eq_prob_princc}:
\begin{equation}
\label{eq_dual_func} d(\nu,\mu)=\min_{z \in \rset^n} \lu
(z,\nu,\mu),
\end{equation}
with  the Lagrangian function $\lu(z,\nu,\mu) = f(z)+\langle \nu,
Az-b \rangle +\langle \mu, Cz-c \rangle$.  For simplicity of the
exposition we introduce further the following notations:
\begin{equation}
\label{eq_notations} G=\left[
\begin{array}{c}
  A \\
  C \\
\end{array}
\right]~\text{and} ~ g=\left[
\begin{array}{c}
  b \\
  c \\
\end{array}
\right].
\end{equation}
Since $f_i$ are strongly convex functions, then $f$ is also strongly convex w.r.t. Euclidian norm $\| \cdot\|$ with convexity parameter $\sigma_{\mathrm{f}}=\min \limits_{i=1,\dots,M} \sigma_i$. Further,  the dual function $d$ is differentiable and its gradient is given by the following expression \cite{NecNed:13}:
\begin{equation*}
\nabla d(\nu,\mu)=Gz(\nu,\mu)-g,
\end{equation*}
where $z(\nu,\mu)$ denotes the unique optimal solution of the
inner problem  \eqref{eq_dual_func}, i.e.:
\begin{equation}
\label{eq_inner_sol} z(\nu,\mu)=\arg\min_{z \in \rset^n}
\lu(z,\nu,\mu).
\end{equation}
Moreover, the gradient $\nabla d$ of the dual function is Lipschitz
continuous w.r.t. Euclideand  norm $\|\cdot\|$, with constant \cite{NecNed:13}:
\begin{equation*}
L_{\text{d}}=\frac{\left\|G\right\|^2}{\sigma_{\mathrm{f}}}.
\end{equation*}
 If
we denote by $\nu_{\Ncal_i}=\left[\nu_j\right]_{j \in \Ncal_i}$ and
by $\mu_{\Ncal_i}=\left[\mu_j\right]_{j \in \Ncal_i}$ we can observe
that the dual function can be written in the following separable
form:
\begin{equation*}
d(\nu,\mu)=\sum_{i=1}^M d_i(\nu_{\Ncal_i},\mu_{\Ncal_i})-\langle
\nu,b\rangle-\langle \mu, c\rangle,
\end{equation*}
with
\begin{align}
\label{eq_separable_dual}
d_i(\nu_{\Ncal_i},\mu_{\Ncal_i})&=\min_{z_i \in \rset^{n_i}}
f_i(z_i)+ \langle
\nu, A_{i} z_i\rangle+\langle \mu, C_{i} z_i\rangle\\
&=\min_{z_i \in \rset^{n_i}} f_i(z_i)+ \sum_{j \in \Ncal_i}
\left\langle A_{ji}^T \nu_j+C_{ji}^T\mu_j,z_i\right\rangle.
\nonumber
\end{align}
In these settings, we have that the gradient $\nabla d_i$ is given
by:
\begin{equation*}
\nabla d_i(\nu_{\Ncal_i},\mu_{\Ncal_i})=\left[\begin{array}{c}
  \left[A_{ji}\right]_{j \in
\Ncal_i} \\
  \left[C_{ji}\right]_{j \in
\Ncal_i} \\
\end{array}\right] z_i(\nu_{\Ncal_i},\mu_{\Ncal_i}),
\end{equation*}
where $z_i(\nu_{\Ncal_i},\mu_{\Ncal_i})$ denotes the unique
optimal solution in \eqref{eq_separable_dual}. Note that $\nabla
d_i$ is Lipschitz continuous  w.r.t. Euclidean norm $\| \cdot \|$,  with constant:
\begin{equation}
L_{\text{d}_i}=\frac{\left\|\left[
\begin{array}{c}
  \left[A_{ji}\right]_{j \in
\Ncal_i} \\
   \left[C_{ji}\right]_{j \in
\Ncal_i} \\
\end{array}
\right]\right\|^2}{\sigma_i}.
\end{equation}
For simplicity of the exposition we will consider further the
notations:
\begin{equation*}
\lambda=\left[\nu^T \mu^T
\right]^T~\text{and}~\lambda_j=\left[\nu_j^T ~\mu_j^T \right]^T
\quad \forall j \in V_2,
\end{equation*}
and we will also denote the effective  domain of the dual function
by $\mathbb{D}=\rset^p \times \rset^q_+$. The following result,
which is a distributed version of descent lemma is central in our
derivations of distributed algorithms and in our proofs  for the
convergence rate for them. Note that a similar result for the case
of inequality constraints can be also found in \cite{BecNed:13}.

\begin{lemma}
\label{descent_lemma} Let Assumption \ref{ass_strong} $(a)$ hold.
Then, the following inequality is valid:
\begin{equation}
\label{eq_descent} d(\lambda) \geq d(\bar{\lambda})+\left\langle
\nabla d(\bar{\lambda}), \lambda- \bar{\lambda}\right\rangle
-\frac{1}{2}\|\lambda-\bar{\lambda}\|_W^2 \quad \forall \lambda,
\bar{\lambda} \in \mathbb{D},
\end{equation}
where the matrix $W=\text{diag}(W_{\nu},W_{\mu})$ with
$W_{\nu}=\text{diag}\left(\sum_{i \in \bar{\Ncal}_j
}L_{\text{d}_i}I_{p_j}; j \in V_2\right)$ and
$W_{\mu}=\text{diag}\left(\sum_{i \in \bar{\Ncal}_j
}L_{\text{d}_i}I_{q_j}; j \in V_2\right)$.
\end{lemma}
\proof Let us first denote by
$\lambda_{\Ncal_i}=\left[\nu_{\Ncal_i}^T~
\mu_{\Ncal_i}^T\right]^T$. Using now the continuous Lipschitz
gradient property of $d_i$ we can write for each $i=1,\dots,M$:
\begin{align*}
d_i(\lambda_{\Ncal_i})\geq
d_i(\bar{\lambda}_{\Ncal_i})+\left\langle \nabla
d_i(\bar{\lambda}_{\Ncal_i}),\lambda_{\Ncal_i}
-\bar{\lambda}_{\Ncal_i}\right\rangle-\frac{L_{\text{d}_i}}{2}\|\lambda_{\Ncal_i}-\bar{\lambda}_{\Ncal_i}\|^2.
\end{align*}
Summing up these inequalities for all $i=1,\dots,M$ and adding
$\langle\lambda, \left[b^T c^T\right]^T\rangle$ to both sides of
the obtained inequality we obtain:
\begin{align}
d(\lambda) \geq d(\bar{\lambda}) + \left\langle \nabla
d(\bar{\lambda}), \lambda - \bar{\lambda} \right \rangle
-\sum_{i=1}^M
\frac{L_{\text{d}_i}}{2}\|\lambda_{\Ncal_i}-\bar{\lambda}_{\Ncal_i}\|^2.
\end{align}

Using now the definitions of $\lambda_{\Ncal_i}$ and $W$ we can
write:
\begin{align*}
\sum_{i=1}^M
\frac{L_{\text{d}_i}}{2}\|\lambda_{\Ncal_i}-\bar{\lambda}_{\Ncal_i}\|^2=\sum_{j=1}^{\bar{M}}
\sum_{i \in \bar{\Ncal}_j}
\frac{L_{\text{d}_i}}{2}\|\lambda_j-\bar{\lambda}_j\|^2=\frac{1}{2}\|\lambda-\bar{\lambda}\|_W.
\end{align*}
Introducing this result into the previous inequality we conclude
the statement. \qed

The following result, which is an extension of Lemma 2.2 in
\cite{BecNed:13} to the case when both equality and inequality
constraints are present, will be useful for characterizing the
distance between a primal estimate and the primal optimal solution
$z^*$ of our optimization problem \eqref{eq_prob_princc}.

\begin{lemma}
\label{lemma_dist_optim} Let Assumption \ref{ass_strong} hold.
Then, the following relation is valid:
\begin{equation}
\label{eq_lemma_dist_optim}
\frac{\sigma_{\mathrm{f}}}{2}\|z(\lambda)-z^*\|^2 \leq f^*
-d(\lambda) \quad \forall \lambda \in \mathbb{D},
\end{equation}
where $z(\lambda)= \arg\min \limits_{z \in \rset^n}
\lu(z,\lambda)$.
\end{lemma}
\proof Since $f$ is $\sigma_{\mathrm{f}}$-strongly convex it
follows that $\lu(z,\lambda)$ is also
$\sigma_{\mathrm{f}}$-strongly convex in the variable $z$ which
together with the definition of $d(\lambda)=f(z(\lambda)) + \langle \lambda,
G z(\lambda) - g \rangle$ and $\nabla d(\lambda)=Gz(\lambda)-g$ and the fact
that $\langle \lambda, \nabla d(\lambda^*) \rangle \leq 0$ gives:
\begin{align*}
\frac{\sigma_{\mathrm{f}}}{2}\|z(\lambda)\!-\!z^*\|&\leq\!
\lu(z^*,\lambda)\!-\!\lu(z(\lambda),\lambda)\!\\
&=\! f(z^*)+\langle \lambda, \nabla
d(\lambda^*)\rangle\!-\!f(z(\lambda))\!-\!\langle \lambda, \nabla
d(\lambda)\rangle\leq f^*-d(\lambda).
\end{align*}\qed

We  denote by $\Lambda^*$ the set of optimal solutions of dual problem
\eqref{eq_dual_prob}. According to Gauvin's theorem \cite{Gau:77},
if Assumption \ref{ass_strong} holds for our original problem
\eqref{eq_prob_princc}, then  $\Lambda^*$ is nonempty and bounded. Since
 the set of optimal Lagrange multipliers is bounded, for any  $\lambda^0
\in \rset^{p+q}$ we can define the following finite quantity:
\begin{equation}
\label{eq_multipleirs_bounded} \mathcal{R}(\lambda^0)=\max
\limits_{\lambda^* \in \Lambda^*} \|\lambda^*-\lambda^0\|_W.
\end{equation}

In this paper  we propose  different distributed dual first order methods for which
 we are interested  in deriving  estimates for both dual and
primal suboptimality and also for primal feasibility violation,
i.e. finding a primal-dual pair $\left(\hat{z},\hat{\lambda}\right)$
such that:
\begin{align}
\label{condition_outer} & ~~~~~ \|\left[G \hat{z}
-g\right]_{\mathbb{D}}\|_{W^{-1}}
 \leq \mathcal{O} (\epsilon), \;\;\;  \| \hat{z}- z^* \|^2 \leq
\mathcal{O} (\epsilon), \\
&-\mathcal{O}(\epsilon)\!\leq \!f(\hat{z}) -\! f^* \!\!\leq\!
\mathcal{O} (\epsilon) ~\text{and} ~f^*-d(\hat{\lambda}) \leq
\mathcal{O} (\epsilon), \nonumber
\end{align}
where $\epsilon$ is a given accuracy.


\subsection{Motivation}
\label{subsec_motivation}

Many engineering applications from networks can be posed as linearly
constrained separable convex optimization problems of type
\eqref{eq_prob_princc}. We will discuss further three such applications,
namely network utility maximization (NUM)
problem, optimal power flow (DC-OPF) problem for a power system and distributed model predictive control (MPC) problem
for networked systems.

\subsubsection{Network utility optimization}
\label{subsec_num}

We consider a network characterized by a bipartite graph
$\mathcal{G}=\left(V_1, V_2,E\right)$,  with $V_1=\left\{1,\dots,M\right\}$
a set of sources, $V_2=\left\{1,\dots,\bar{M}\right\}$ a set of
capacitated links, each link $j$ having capacity $\bar{c}_j
> 0$, and $E$ its incidence matrix. In these settings, $\Ncal_i$ represents the set of links which are
used by the source $i$, while $\bar{\Ncal}_j$ is the set of
sources which share the link $j$. Also, we attach to each source
$i$ a strongly convex decreasing utility function $f_i(z_i)$,
where $z_i \in \rset$ denotes the rate at which the source sends
its data. In these settings, the goal of the network utility
problem is to find the optimal rates at which the total utility
function is minimized. Introducing the notation $z=\left[z_1
\cdots z_M\right]^T$, the network utility maximization problem can
be posed as the following convex optimization problem:
\begin{align}
\label{eq_utility} f^*=&\min_{z_i \in \rset}
f(z)~~~\left(=\sum_{i=1}^M
f_i(z_i)\right)\\
&\text{s.t.:}~~\sum_{i \in \bar{\Ncal}_j} z_i \leq \bar{c}_j ~~
\forall j \in V_2,~~~ z_i \in Z_i =[0, \; R_i] ~~ \forall i \in V_1.
\nonumber
\end{align}
By stacking together all the local and coupling constraints, we
can observe that problem \eqref{eq_utility} can be written in the
form of problem \eqref{eq_prob_princc}, where the equality
constraints are absent.
Well known applications are the NUM problem \cite{BecNed:13} and
dynamic network utility maximization (DNUM) with end-to-end delays \cite{TriZym:08}.\\

\subsubsection{DC Optimal power flow}
\label{subsec_opf}

Let us discuss the active optimal power flow
(DC-OPF) problem for a power system \cite{ZimMur:11}. We consider a
power system whose structure is characterized by a directed
bipartite graph $\mathcal{G}=(V_1,V_2,E)$, where $V_1=\{i~|~i=1,\dots,M\}$
denotes the set of buses, $V_2=\{l=(i,j)~|~i,j \in V_1,
~l=1,\dots,\bar{M}\} \subseteq V_1 \times V_1$ represent the sets of
transmission lines (branches) between two buses and the matrix $E$
denotes its incidence matrix. In these settings we have:
\begin{align*}
 \Ncal_i & =\{l \in V_2~|~ E_{li} \neq 0\} =\{l \in V_2~|~ \exists j \in V_1 ~\text{s.t.}~(i,j)\vee(j,i)=l\}
\end{align*}
which  denotes the set of all transmission lines from or to bus $i$
and
\begin{align*}
\bar{\Ncal}_l & =\{i \in V_1~|~ E_{li} \neq 0\} =\{i,j \in V_1~|~(i,j)\vee (j,i)=l\}
\end{align*}
which  denotes the set comprised of buses $i$ and $j$ which define
the branch $l$. We also introduce:
\begin{align*}
\Scal_i  =\bigcup_{l \in V_2}\{j \in V_1~|~ E_{li}\neq 0
\wedge E_{lj} \neq 0\} =\{j \in V_1 ~|~ \exists l \in V_2~\text{s.t.}~(i,j)\vee (j,i)=l\}
\end{align*}
which denotes the sets of all buses directly linked with bus $i$. It
is straightforward to notice that the set $\Scal_i$ can be
obtained from the sets $\Ncal_i$ and $\bar{\Ncal}_l$.

We define further the diagonal matrix $R \in \rset^{\bar{M}\times
\bar{M}}$, whose diagonal elements $R_{ll}$ represent the reactance
of the $l$th transmission line between two busses $i$ and $j \in
V_1$. For each bus $i$ we denote by \vspace*{-0.2cm}\[ \theta_i \in \Theta_i=
[\underline{\theta}_i,\overline{\theta}_i] \]  the phase angle of
the voltage and by  \[ P_i^g \in \mathcal{P}_i =
[\underline{P}_i^g,\overline{P}_i^g] \] the generated power if the
bus $i$ is directly connected to a generator. Under this model, the
active power flow from a bus $i$ to a bus $j$ is given by:
\vspace{-0.2cm}
\begin{equation}
\label{eq_active_flow}
F_{l}=\frac{1}{R_{ll}}\left(\theta_i-\theta_j\right),
\end{equation}
where $l=(i,j)$ and we recall that $R_{ll}$ represent the reactance
of the transmission line connecting buses $i$ and $j$. We impose
lower and upper line flows limits
$\underline{F}=\left[\underline{F}_1 \cdots
\underline{F}_{\bar{M}}\right]^T$ and
$\overline{F}=\left[\overline{F}_1 \cdots
\overline{F}_{\bar{M}}\right]^T$, respectively. We also assume that
each bus $i$ is characterized by a local load $P_i^d$ and we denote
by $P^d=\left[P_1^d \cdots P_M^d\right]^T$ the overall vector of
loads. We introduce further the notations:
\vspace*{-0.1cm}
\begin{equation*}
\theta=\left[\theta_1 \cdots
\theta_M\right]^T~\text{and}~P^g=\left[P_1^g\cdots
P_{M_g}^g\right]^T,
\end{equation*}
where $M_g$ denotes the number of generators. We also define the
matrix $A^g \in [0,1]^{M \times M_g}$ having $A^g_{ij}=1$ if $P_j^g$
is directly linked with the bus $i$ and the rest of its entry equal
to zero. Note that if we consider that each bus $i$ is directly
coupled with a generator unit $A^g=I_{M}$. Using these notations,
the DC nodal power balance can be written in the following form
\cite{ZimMur:11}:
\begin{equation}
\label{eq_nodal} E^TRE \theta =A^gP^g- P^d,
\end{equation}
where the matrix $E^TRE$ denotes the weighted Laplacian and its entries have the
following expressions:
\begin{equation*}
 [E^TRE]_{ij}=\left\{ \begin{array}{c}
                       \sum_{s \in \Scal_i} R_{ll}, ~l=(i,s)\vee(s,i)~~\text{if}~i=j \\
                       -R_{ll}, ~l=(i,j)\vee(j,i)~~~~~~~~\text{if}~i \neq j\\
                       0~~~~~~~~~~~~~~~~~~~~~~~~~~~~~~~~\text{otherwise}.
                     \end{array} \right.
\end{equation*}
We can observe that the structure of the Laplacian matrix is given
by the structure of the incidence matrix $E$ through the sets
$\bar{\Scal}_i$, which, at its turn depend on the sets $\Scal_i$ and
$\Ncal_l$ for all $i \in V_1$ and $l \in V_2$. Using further the relation between the the power flow and the phase
angle of the voltages, we can write the lower and upper limits
imposed on the line flows in the following matrix form:
\begin{equation}
\label{eq_limit_flow} \underline{F} \leq RE \theta \leq
\overline{F}.
\end{equation}
We also define reference values $\theta_i^{\text{ref}}$ for the
phase angle of the voltage of each bus and $P_i^{g,\text{ref}}$ for
the generated powers of each generator. Further, for each bus $i$ we
define a local decision variable $z_i$ as follows:
\begin{equation*}
z_i=\left\{\begin{array}{c}
            \left[\begin{array}{c}
              \theta_i \\
              P_i^g
            \end{array}\right]
            ~\text{if the bus} \; i \;  \text{is connected to a generator}\\
             \!\!\!\!\!\!\!\!\!\!\!\!\!\!\!\!\!\!\!\!\!\!\!\!\!\!\!\!\!\!\!\!\theta_i \qquad \qquad\ \qquad \text{otherwise}
          \end{array}\right.
\end{equation*}
and the corresponding reference values $z_i^{\text{ref}}$.

In comparison with the approach made in \cite{ZimMur:11}, where the
authors consider the lower and upper limits of the form
$\theta_i^{\text{ref}} \leq \theta_i \leq \theta_i^{\text{ref}}$, in
our approach we do not impose such constraints but use instead a
weighted quadratic cost, which, depending on the value of the
parameter $q_i$, requires  the solution to be close to the reference
value $\theta_i^{\text{ref}}$. The main motivation behind this
approach consist in the fact that constraints of this form usually
induce numerical problems due to the fact that the optimization
problem which has to be solved is badly conditioned (for example,
the Slater constraint qualification does not hold in this case).
Therefore, for each bus $i$ directly connected to a generator unit
we impose a local cost of the form:
\begin{equation}
\label{eq_cost_local}
f_i(z_i)=0.5\|z_i-z_i^{\text{ref}}\|_{Q_i}^2-\gamma_i \log(\beta_i+
P_i^g),
\end{equation}
where the diagonal matrix $Q_i=\left[
             \begin{array}{cc}
               q_i & 0 \\
               0 & p_i \\
             \end{array}
           \right]
\in \rset^{2 \times 2}$ and the positive scalar $\gamma_i$ are used
in order to weight the local cost. Also, the positive scalar
$\beta_i$ is used to avoid numerical instability when $P_i^g$ is
closed to $0$. Also, in comparison with the existing approaches for
(DC-OPF) problems we add to the classic quadratic term a weighted
logarithmic term, which is used in many resource allocation problem
(see e.g. \cite{XiaBoy:06}) in order to reduce the absolute risk aversion. The
logarithmic utility function also exhibit diminishing returns with the rate of resources, in our case
the generated power, that is, as rate increases the incremental utility grows by smaller amounts.
For the buses that are not connected to a generator unit we impose a simple quadratic local cost of the
form:
\begin{equation}
\label{eq_cost_local_simple}
f_i(z_i)=0.5q_i\left(z_i-\theta_i^{\text{ref}}\right)^2,
\end{equation}
where in this case $q_i$ is a positive scalar. Note that for these
choices the local costs $f_i$ are strongly convex functions for both
cases. In conclusion, the (DC-OPF) problem can be cast as the following large-scale
separable convex optimization problem:
\begin{align}
\label{eq_prob_opf} f^* = &\min_{\theta_{i} \in \Theta_i, P^g_i \in \mathcal{P}_i}
\sum_{i_1} f_{i_1}(\theta_{i_1})+\sum_{i_2} f_{i_2}(\theta_{i_2},P^g_{i_2})\\
&\text{s.t.:} ~~E^TRE \theta-A^{g}P^{g}=-P^d, ~~\underline{F} \leq RE \theta \leq \overline{F}. \nonumber
\end{align}

\subsubsection{Distributed MPC for networked systems}
\label{subsec_mpc}

We consider a discrete-time networked system, modelled  by a graph
$\mathcal{G}=\left(V,E\right)$, for which the set
$V=\left\{1,\dots,M\right\}$ represents the subsystems and the
adjacency matrix $E$ indicates the dynamic couplings between these
subsystems. The dynamics of the subsystems can be defined by the
following linear state equations \cite{NecSuy:08}:
\begin{equation}
\label{mod3} x_i(t+1) =\sum _{j \in \mathcal {N}_{i} }
\bar{A}_{ij} x_j(t)+
 \bar{B}_{ij} u_j(t)  \qquad \forall i \in V,
\end{equation}
where $x_i(t) \in \mathbb{R}^{n_{x_i}}$ and $u_i(t) \in
\mathbb{R}^{n_{u_i}}$ represent the state and  the input of $i$th
subsystem at time $t$, $\bar{A}_{ij} \in \mathbb{R}^{n_{x_i}
\times n_{x_j}}$ and $\bar{B}_{ij} \in \mathbb{R}^{n_{x_i} \times
n_{u_j}}$. Note that in these settings $\mathcal {N}_{i}$ denotes
the set of subsystems, including $i$, whose dynamics directly
affect the dynamics of subsystem $i$ and $\bar{\Ncal}_i$
represents the set of subsystem, including $i$, whose dynamics are
affected by the dynamics of subsystem $i$. We also impose local
state and input constraints: \[ x_i(t) \in X_i, \quad   u_i(t) \in
U_i \qquad \forall i \in V, \;\; t \geq 0,
\]
where $X_i \subseteq \rset^{n_{x_i}}$ and $U_i \subseteq
\rset^{n_{u_i}}$ are polyhedral sets. For a prediction horizon of
length $N$, we consider strongly convex stage and final costs for
each subsystem $i$:
\[ \sum_{t=0}^{N-1} \ell_i(x_i(t),u_i(t)) + \ell_i^{\mathrm{f}}(x_i(N)),
\]
where the final costs $\ell_i^{\mathrm{f}}$ are chosen  such that
the control scheme ensures stability
\cite{ScoMay:99,NecNed:13,NecSuy:08}. The centralized MPC problem
for the networked system \eqref{mod3}, for a given initial state
$x=\left[x_1^{T}\cdots x_M^{T}\right]^T$ can be posed as the
following convex optimization problem:
\begin{align}
 &  \min _{x_i(t),u_i(t)}  \sum_{i=1}^M \sum_{t=0}^{N-1}
\ell_i(x_i(t),u_i(t)) + \ell_i^{\mathrm{f}}(x_i(N))  \nonumber \\
& \text{s.t.:} \;\; x_i(t+1) =
\sum _{j \in \mathcal {N}^{i} } \bar{A}_{ij}  x_j (t) + \bar{B}_{ij} u_j(t), \; x_i(0) =x_i,  \label{mod5} \\
& \;\;\;\;\;\;\;\;  x_i(t) \in X_i, \; u_i(t) \in U_i, \; x_i(N)
\in X_i^{\mathrm{f}}   \quad  \forall i \in V, ~ \forall t,
\nonumber
\end{align}
where $X_i^{\mathrm{f}}$ are terminal sets  chosen under some
appropriate conditions to ensure stability of the MPC scheme (see
e.g. \cite{ScoMay:99,NecNed:13,NecSuy:08}). For the state and
input trajectory of subsystem $i$ and the overall state and input
trajectory we use the notations:
\begin{align*}
z_i& = \left[u_{i}(0)^T x_i(1)^T \cdots u_{i}(N-1)^T
x_i(N)^T \right]^T,\\
z&=\left[z_{1}^T \cdots z_{M} ^T\right]^T,
\end{align*}
and for the total local cost over the prediction horizon and local
constraints of each subsystem we introduce:
\begin{align*}
&f_i(z_i)=\sum_{t=0}^{N-1} \ell_i(x_i(t),u_i(t)) +
\ell_i^{\mathrm{f}}(x_i(N)),\\
&Z_i=\left(\prod_{i=1}^{N-1} U_i \times X_i\right) \times U_i
\times X_i^\mathrm{f}.
\end{align*}

In these settings, the optimization problem \eqref{mod5} can be
written equivalently as the structured optimization problem
\eqref{eq_prob_princc} where $n_i=N(n_{u_i}+n_{x_i})$, the
equality constraints $Az=b$ are obtained by stacking all the
dynamics \eqref{mod3} together, while the inequality constraints
$Cz \leq c$ are obtained by writing the local constraints $z_i \in
Z_i$ in compact form. Note also that for the matrix $A$, each
block $A_{ij}=0$ whenever $E_{ij}=0$ and $C$ is a block diagonal
matrix.

In the following sections we will propose and analyze dual
distributed  (fast) gradient methods for solving the dual problem
\eqref{eq_dual_prob} which exploit the separability of the the
dual function and allow us to recover a suboptimal and nearly
feasible solution for our original problem \eqref{eq_prob_princc}.


\section{Distributed dual fast gradient algorithm (\textbf{DFG})}
\label{sec_ddfg}

In this section we propose a distributed dual fast (also called
accelerated) gradient scheme (\textbf{DFG}) for solving the dual
problem \eqref{eq_dual_prob}. A similar algorithm was proposed by
Nesterov in \cite{Nes:05} and applied further in \cite{NecSuy:08}
for solving dual problems. A similar version of the algorithm was
also proposed in \cite{NecNed:13} for the case when the dual updates
use inexact information and the step size is a fixed scalar. The
scheme defines two sequences $\left({\hat
\lambda}^k,\lambda^k\right)_{k\geq 0}$ for the dual variables:
\begin{center}
\framebox{
\parbox{8.3cm}{
\begin{center}
\textbf{ Algorithm {\bf (DFG)}}
\end{center}
{Initialization: $\lambda^0=0$. For $k\geq 0$ compute:}
\begin{enumerate}
\item{$z^k = \arg \min\limits_{z \in \rset^n}
\mathcal{L}(z,\lambda^k)$}

\item ${\hat \lambda}^k=\left[\lambda^k+W^{-1}\nabla
d(\lambda^k)\right]_{\mathbb{D}}$

\item $\lambda^{k+1}=\frac{k+1}{k+3}\!{\hat
\lambda}^k\!+\!\frac{2}{k+3}\!\left[W^{-1}\!\sum_{s=0}^k
\!\frac{s+1}{2} {\bar \nabla} d(\lambda^s)\right]_{\mathbb{D}}$.
\end{enumerate}
}}
\end{center}

For simplicity of the exposition we restrict our analysis to the
case $\lambda^0=0$. Note that the behavior of the Algorithm
(\textbf{DFG}) remain unchanged if one use any initialization
$\lambda^0 \in \mathbb{D}$ (see e.g. \cite{NecNed:13}).  We can
also observe that step $1$ of the algorithm requires an exact
solution of the inner optimization problem. In many practical
applications such a solution is hard to be computed. Instead, one
can compute an approximate solution, i.e. $\bar{z}^k \approx \arg
\min_{z \in \rset^n} \mathcal{L}(z,\lambda^k)$,  which satisfies a
certain inner accuracy (see  \cite{NecNed:13} for a detailed
discussion). The main difference between our Algorithm
(\textbf{DFG}) and the algorithms proposed in
\cite{Nes:05,NecSuy:08,NecNed:13} consists in the way we update
the sequence $\lambda^k$. Instead of using a classical projected
gradient step with a scalar step size as in
\cite{Nes:05,NecSuy:08,NecNed:13}, we update $\lambda^k$ using a
projected weighted gradient step which allows us to obtain a
distributed scheme. Further, we will analyze the convergence
properties of  Algorithm (\textbf{DFG}).

\subsection{Sublinear convergence using an average primal sequence}
\label{subsec_convergence_dfwg}

As we have stated before, in this section we are interested in
characterizing the dual suboptimality and also the primal
suboptimality and feasibility violation for Algorithm
(\textbf{DFG}). Using Lemma \ref{descent_lemma} instead of the
classical descent lemma we can obtain from Theorem 2 in
\cite{Nes:05} the following inequality, which will help us to
establish the convergence properties of Algorithm (\textbf{DFG}):
\begin{align}
\label{eq_rec} &\frac{(k+1)(k+2)}{4} d({\hat \lambda}^k) \nonumber\\
&~~\geq \max_{\lambda \in \mathbb{D}}-\frac{1}{2}
\|\lambda\|_{W}^2 \!+\!\sum_{s=0}^k
\!\frac{s\!+\!1}{2}\!\left[d(\lambda^s)\!+\!\langle  \nabla
d(\lambda^s),\lambda\!-\!\lambda^s \rangle\right] \quad \forall
\lambda \in \mathbb{D}.
\end{align}

The following theorem provides an estimate on the dual
suboptimality for Algorithm ({\bf DFG}):
\begin{theorem}
\label{theorem_dual_optim} Let Assumption \ref{ass_strong} hold
and the sequences $\left(z^k,{\hat
\lambda}^k,\lambda^k\right)_{k\geq 0}$ be generated by algorithm
({\bf DFG}). Then, an estimate on dual suboptimality for
\eqref{eq_dual_prob} is given by:
\begin{equation}
\label{bound_dual_optim} f^*-d({\hat \lambda}^k)\leq
\frac{2\mathcal{R}^2}{(k+1)^2},
\end{equation}
where $\mathcal{R}=\mathcal{R}(0) =\max \limits_{\lambda^* \in
\Lambda^*} \|\lambda^*\|_W$ according to
\eqref{eq_multipleirs_bounded}.
\end{theorem}

\proof Using the concavity of $d$ and $\lambda=\lambda^*$ in
\eqref{eq_rec} we get:
\begin{align*}
\frac{(k+1)(k+2)}{4} d({\hat \lambda}^k) \geq &- \frac{1}{2}
\|\lambda^*\|^2_W+\sum_{s=0}^k\frac{s+1}{2}d(\lambda^*).
\end{align*}
Dividing now both sides  by $\frac{(k+1)(k+2)}{4}$, rearranging
the terms and taking into account that $d(\lambda^*)=f^*$,
$(k+1)^2 \leq (k+1)(k+2)$ and the definition of $\mathcal{R}$ we
obtain \eqref{bound_dual_optim}. \qed

We define further the following average sequence for the primal
variables:
\begin{equation}
\label{primal_point} {\hat
z}^k=\sum_{s=0}^k\frac{2(s+1)}{(k+1)(k+2)}z^s.
\end{equation}

The next result gives an estimate on primal feasibility violation.
\begin{theorem}
\label{theorem_primal_fesa} Under the assumptions of Theorem
\ref{theorem_dual_optim} and ${\hat z}^k$ generated by
\eqref{primal_point}, an estimate on primal feasibility violation
for original problem \eqref{eq_prob_princc} is given by:
\begin{equation}
\label{ineq_bound_fesa_fg} \left\|\left[
\begin{array}{c}
  A \hat{z}^k-b \\
  ~~\left[C \hat{z}^k-c\right]_{\rset^q_+} \\
\end{array}%
\right]\right\|_{W^{-1}} \leq \frac{8\mathcal{R}}{(k+1)^2}.
\end{equation}
\end{theorem}

\proof Using \eqref{eq_rec}, the convexity of $f$ and the
definitions of $d$ and $\nabla d$ , we can write for any $\lambda
\in \mathbb{D}$:
\begin{align}
\label{ineq_fesa0} \max_{\lambda \in \mathbb{D}
}-\frac{2}{(k+1)^2}\|\lambda\|^2_W+\langle \lambda, G \hat{z}^k-g
\rangle \leq d({\hat \lambda}^k) - f({\hat z}^k).
\end{align}
For the second term of the right-hand side we have:
\begin{align}
\label{ineq_fesa1}  d({\hat \lambda}^k)-f({ \hat z}^k)&\leq
d(\lambda^*)-f({ \hat z}^k) = \min_{z \in \rset^n} f(z) +\langle
\lambda^*, Gz-g \rangle
-f({\hat z}^k)\\
&\leq \!f({ \hat z}^k)\!+\!\!\langle \lambda^*\!,\! G
\hat{z}^k-g\rangle \!-\!\!f({ \hat z}^k)\!=\!\langle
\lambda^*\!,\! G \hat{z}^k-g \rangle \nonumber \leq
\!\langle\lambda^*\!,\!\left[G
\hat{z}^k-g\right]_{\mathbb{D}}\rangle,
\end{align}
where in the last inequality we use that $\lambda^* \in
\mathbb{D}$. By evaluating the maximum in the left-hand side term
in \eqref{ineq_fesa0} and taking into account that $\langle [G
\hat{z}^k-g ]_{\mathbb{D}},G \hat{z}^k-g -[G \hat{z}^k-g
]_{\mathbb{D}}\rangle=0$ we obtain the following relation:
\begin{align}
\label{ineq_fesa2} \max_{\lambda \in
\mathbb{D}}-&\frac{2}{(k+1)^2}\|\lambda\|^2_W+\langle \lambda,
G\hat{z}^k-g \rangle = \frac{(k+1)^2}{8}\|[G \hat{z}^k-g
]_{\mathbb{D}}\|^2_{W^{-1}}.
\end{align}
Combining now \eqref{ineq_fesa1} and \eqref{ineq_fesa2} with
\eqref{ineq_fesa0}, using the Cauchy-Schwartz inequality and
introducing the notation $\alpha=\|\left[G
\hat{z}^k-g\right]_{\mathbb{D}}\|_{W^{-1}}$, we obtain the
following second order inequality in $\alpha$:
\begin{equation*}
\frac{(k+1)^2}{8}\alpha^2-\|\lambda^*\|_{W}\alpha \leq 0,
\end{equation*}
from which, using the definitions of $G \hat{z}^k-g$, $\mathbb{D}$
and $\mathcal{R}$ we get \eqref{ineq_bound_fesa_fg}. \qed

\begin{theorem}
\label{theorem_primal_optim} Assume that the conditions in Theorem
\ref{theorem_primal_fesa} are satisfied and let ${\hat z}^k$ be
given by \eqref{primal_point}. Then, the following estimate on
primal suboptimality  for  problem \eqref{eq_prob_princc} can be
derived:
\begin{equation}
\label{bound_primal_optim} -\frac{8 \mathcal{R}^2}{(k+1)^2}\leq\!
f(\hat{z}^k) -\! f^* \!\!\leq 0.
\end{equation}

Moreover, the sequence $\hat{z}^k$ converges to the unique optimal
solution $z^*$ of \eqref{eq_prob_princc} with the the following
rate:
\begin{equation}
\label{bound_primal_solution} \|\hat{z}^k-z^*\| \leq \frac{4
\mathcal{R}}{\sqrt{\sigma_{\mathrm{f}}}(k+1)}.
\end{equation}

\end{theorem}

\proof The right-hand side inequality in
\eqref{bound_primal_optim} follows from evaluating
\eqref{ineq_fesa0} in $\lambda=0$ and taking into account that
$d({\hat \lambda}^k) \leq d(\lambda^*)=f^*$.

\noindent In order to prove the left-hand side inequality of
\eqref{bound_primal_optim} we can write:
\begin{align}
\label{ineq_aux_optim} f^*=d(\lambda^*)&=\min_{z\in \rset^n}
f(z)+\langle
\lambda^*,Gz-g\rangle \\
&\leq f({ \hat z}^k)+\langle \lambda^*,G \hat{z}^k-g \rangle \nonumber \\
&\leq f({\hat z}^k)+
\langle \lambda^*,\left[G \hat{z}^k-g\right]_{\mathbb{D}} \rangle \nonumber \\
&\leq f({ \hat z}^k) + \|\lambda^*\|_{W}\left\|\left[G
\hat{z}^k-g\right]_{\mathbb{D}}\right\|_{W^{-1}}, \nonumber
\end{align}
where the second inequality follows from the fact that $\lambda^*
\in \mathbb{D}$ and the last one from Cauchy-Schwartz inequality.
Using now \eqref{ineq_bound_fesa_fg} we obtain the result.

Further, since $f$ is $\sigma_{\mathrm{f}}$-strongly convex, we
have also that $\lu(z,\lambda)$ is also
$\sigma_{\mathrm{f}}$-strongly convex for all $\lambda \in
\mathbb{D}$. Thus, taking $\lambda=\lambda^*$ and noting that
$z^*=\arg\min_{z \in \rset^n} \lu(z,\lambda^*)$ we have:
\begin{align*}
\frac{\sigma_{\mathrm{f}}}{2}\|\hat{z}^k-z^*\|^2 &\leq
\lu(\hat{z}^k,\lambda^*)-\lu(z^*,\lambda^*)\\
&=f(\hat{z}^k)+\langle \lambda^*, G \hat{z}^k-g \rangle -f^*
\overset{\eqref{bound_primal_optim}}{\leq} \langle \lambda^*, G
\hat{z}^k-g \rangle\\
&\leq \langle \lambda^*,\left[G \hat{z}^k-g\right]_{\mathbb{D}}
\rangle \leq \|\lambda^*\|_{W}\left\|\left[G
\hat{z}^k-g\right]_{\mathbb{D}}\right\|_{W^{-1}},
\end{align*}
where the last two inequalities follows from the same arguments as
in \eqref{ineq_aux_optim}. Using now \eqref{ineq_bound_fesa_fg}
and the definition of $\mathcal{R}$ we obtain
\eqref{bound_primal_solution}. \qed

\begin{remark}
\label{remark_dfg} $(i)$ If we use for the initialization of the
algorithm any $\lambda^0 \in \mathbb{D}$ the order of the
estimates on primal and dual suboptimality and primal feasibility
violation derived above remain unchanged.\\
$(ii)$ Note that according to Theorem \ref{theorem_primal_optim}
for $\lambda^0=0$ we are always below the optimal value $f^*$. In
the case when we use an initialization $\lambda^0 \neq
0$ we cannot guarantee anymore this property.\\
$(iii)$ From previous theorems we observe that for a given accuracy
$\epsilon$, we  need to perform
$\mathcal{O}\left(\frac{1}{\sqrt{\epsilon}}\right)$ iterations in
order to obtain  a primal suboptimal and near-feasible solution
based on averaging the primal generated sequence. \qed
\end{remark}


\section{Hybrid distributed dual fast gradient algorithm (\textbf{H-DFG})}
\label{sec_hddfg}

Note that for the Algorithm (\textbf{DFG}) the primal sequence
$\left\{\hat{z}^k\right\}_{k \geq 0}$ for which we are able to
recover primal suboptimality and primal feasibility violation is
given by a weighted average of the iterates $\left\{z^k\right\}_{k
\geq 0}$.  However, in simulations we observe also a good
behaviour of the last iterate $z^k$. In this section we propose a
hybrid distributed dual fast gradient algorithm for which we can
ensure estimates for both primal suboptimality and feasibility
violation of the last iterate $z^k$, which  supports our findings
from simulations. The algorithm is characterized by two phases: in
the first phase we perform $k$ steps of Algorithm (\textbf{DFG})
while in the second phase another $k$ steps of a dual weighted
gradient algorithm are performed. A similar hybrid strategy was
also discussed in \cite{Nes:12a,NecNes:11}. We present further the
proposed scheme:
\begin{center}
\framebox{
\parbox{8.3cm}{
\begin{center}
\textbf{ Algorithm {\bf (H-DFG)}}
\end{center}
{Initialization: $\lambda^0=0$. \\
{ \bf Phase 1:} For $j=0,\dots,k$ compute:}
\begin{enumerate}
\item{$z^j = \arg \min\limits_{z \in \rset^n}
\mathcal{L}(z,\lambda^j)$}

\item ${\hat \lambda}^j=\left[\lambda^j+W^{-1}\nabla
d(\lambda^j)\right]_{\mathbb{D}}$

\item $\lambda^{j+1}=\frac{j+1}{j+3}\!{\hat
\lambda}^j\!+\!\frac{2}{j+3}\!\left[W^{-1}\!\sum_{s=0}^j
\!\frac{s+1}{2} {\bar \nabla} d(\lambda^s)\right]_{\mathbb{D}}$.
\end{enumerate}
{{\bf Phase 2:} Set $\lambda^k=\hat{\lambda}^k$. For
$j=k,\dots,2k$ compute:}
\begin{enumerate}
\item{$z^j = \arg \min\limits_{z \in \rset^n}
\mathcal{L}(z,\lambda^j)$}

\item ${\lambda}^{j+1}=\left[\lambda^j+W^{-1}\nabla
d(\lambda^j)\right]_{\mathbb{D}}.$
\end{enumerate}
}}
\end{center}

The following lemma, which is a generalization of a standard
result for gradient methods shows that Phase 2 of Algorithm
(\textbf{H-DFG}) is an ascent method. For completeness we also
give the proof.

\begin{lemma}
\label{lemma_dual_ascent} Let the sequence $\left\{\lambda^j
\right\}_j$ be generated by the Phase 2 of Algorithm
(\textbf{H-DFG}). Then, the value of the dual function increases
at each iteration according to the following relation:
\begin{equation}
\label{eq_ascent_dual} d(\lambda^{j+1}) \geq
d(\lambda^j)+\frac{1}{2}\|\lambda^j-\lambda^{j+1}\|_W^2  \quad  \forall
j=k,\dots,2k.
\end{equation}
\end{lemma}
\proof Let us first notice that the update
${\lambda}^{j+1}=\left[\lambda^j+W^{-1}\nabla
d(\lambda^j)\right]_{\mathbb{D}}$ can be posed as the minimization
of the following quadratic approximation of $d$:
\begin{equation}
\label{eq_projection_dual} \lambda^{j+1}=\arg\max \limits_{\lambda
\in \mathbb{D}} d(\lambda^j)+\langle \nabla d(\lambda^j),
\lambda-\lambda^j \rangle -\frac{1}{2}\|\lambda-\lambda^j\|_W^2.
\end{equation}
From the optimality conditions of problem
\eqref{eq_projection_dual} we obtain:
\begin{equation}
\label{eq_optimality_lambda} \langle \nabla d(\lambda^j),
\lambda^{j+1}-\lambda^j \rangle \geq
\|\lambda^j-\lambda^{j+1}\|_W^2.
\end{equation}
Using now this inequality in Lemma \ref{descent_lemma} with
$\bar{\lambda}=\lambda^j$ and $\lambda=\lambda^{j+1}$ we obtain
the result. \qed

\subsection{Sublinear convergence using the last primal iterate}
\label{subsec_convergence_h-dfg}

We introduce further the following notation:
\begin{equation}
\label{eq_kstar} k^*=\arg\min \limits_{j \in \left[k,2k\right]}
\|\lambda^j-\lambda^{j+1}\|_W^2.
\end{equation}
Note that the quantity $\lambda^j-\lambda^{j+1}$ denotes the
constrained gradient direction (see  \cite{Nes:04}), which
represent an indicator for the suboptimality level of the estimate
$\lambda^j$. We can also observe that $\lambda^j$ is an optimal
solution of \eqref{eq_dual_prob} if and only if
$\lambda^j-\lambda^{j+1}=0$ and thus we want
$\|\lambda^j-\lambda^{j+1}\|_W^2$ to be small. The following
theorem gives an estimate on the dual suboptimality for the
Algorithm (\textbf{H-DFG}):

\begin{theorem}
\label{th_dual_suboptimality_hdfwg} Let Assumption
\ref{ass_strong} hold, the sequences
$\left\{\lambda^j,\hat{\lambda}^j,z^j\right\}_{j \geq 0}$ be
generated by the Algorithm (\textbf{H-DFG}) and $k^*$ be given by
\eqref{eq_kstar}. Then, an estimate for dual suboptimality for
\eqref{eq_dual_prob} is given by:
\begin{equation*}
f^*-d(\lambda^{k^*}) \leq \frac{2\mathcal{R}^2}{(k+1)^2}.
\end{equation*}
\end{theorem}
\proof From Theorem \ref{theorem_dual_optim} and the
initialization in Phase 2 of Algorithm (\textbf{H-DFG}) we have:
\begin{equation*}
\frac{2\|\lambda^*\|_W^2}{(k+1)^2} \geq
f^*-d(\hat{\lambda}^k)=f^*-d(\lambda^k).
\end{equation*}
Using now Lemma \ref{lemma_dual_ascent} we obtain the following
inequalities:
\begin{equation}
\label{eq_decresing} d(\hat{\lambda}^k)=d(\lambda^k) \leq
d(\lambda^{k+1}) \leq \cdots \leq d(\lambda^{2k+1}),
\end{equation}
from which, together with the previous inequality and the fact
that $k^* \in [k,2k]$ we obtain the result. \qed

The following result characterizes the primal feasibility
violation for Algorithm (\textbf{H-DFG}) in the last iterate
$z^{k^*}$.

\begin{theorem}
\label{th_primal_fesability_hdwfg} Under the assumptions of
Theorem \ref{th_dual_suboptimality_hdfwg}, an estimate on primal
feasibility violation for original problem \eqref{eq_prob_princc}
is given by:
\begin{equation}
\label{ineq_bound_fesa_fgh} \left\|\left[
\begin{array}{c}
  A z^{k^*}-b \\
  ~~\left[C z^{k^*}-c\right]_{\rset^q_+} \\
\end{array}%
\right]\right\|_{W^{-1}} \leq
\frac{2\mathcal{R}}{(k+1)\sqrt{(k+1)}}.
\end{equation}
\end{theorem}
\proof Using Theorem \ref{theorem_dual_optim} and Lemma
\ref{lemma_dual_ascent} we can write:
\begin{align*}
\frac{2\|\lambda^*\|_W^2}{(k+1)^2}&\geq
f^*-d(\hat{\lambda}^k)=f^*-d(\lambda^k)\\
&\geq
\!f^*\!\!-\!d(\lambda^{k+1})\!+\!\frac{1}{2}\|\lambda^k\!\!-\!\lambda^{k+1}\|_W^2\!\geq
\cdots\!\geq
f^*\!\!-\!d(\lambda^{2k+1})\!+\!\frac{1}{2}\!\sum_{j=k}^{2k}\|\lambda^j\!\!-\!\lambda^{j+1}\|_W^2\\
&\geq \frac{(k+1)}{2}\|\lambda^{k^*}\!\!-\!\lambda^{k^*+1}\|_W^2,
\end{align*}
where in the last inequality we used \eqref{eq_kstar}. Using the
previous inequality we obtain:
\begin{equation}
\label{eq_bound_grad_direction}
\|\lambda^{k^*}\!\!-\!\lambda^{k^*+1}\|_W^2 \leq
\frac{4\|\lambda^*\|_W^2}{(k+1)^3}.
\end{equation}
Further, we will show that $\left\|\left[\nabla
d(\lambda^{k^*})\right]_{\mathbb{D}}\right\|_{W^{-1}}^2 \leq
\|\lambda^{k^*}\!\!-\!\lambda^{k^*+1}\|_W^2$. We will prove this
inequality componentwise. First, we recall that
$\mathbb{D}=\rset^p \times \rset^q_+$. Thus, for all $i=1,\dots,p$
we have:
\begin{align}
\label{eq_part1} \left|\left[\nabla_i
d(\lambda^{k^*})\right]_{\rset}\right|_{W_{ii}^{-1}}^2\!\!&=\left|\nabla_i
d(\lambda^{k^*})\right|_{W_{ii}^{-1}}^2=\left|\lambda^{k^*}_i-\lambda^{k^*}_i-W_{ii}^{-1}\nabla_i
d(\lambda^{k^*})\right|_{W_{ii}}^2\\
&=\left|\lambda^{k^*}_i-\lambda^{k^*+1}_i\right|_{W_{ii}}^2,
\nonumber
\end{align}
where in the last inequality we used the definition of
$\lambda^{k^*+1}$. We introduce now the following disjoint sets:
$I_{-}=\left\{i \in [p+1,p+q] : \nabla_i d(\lambda^{k^*}) < 0
\right\}$ and $I_{+}=\left\{i \in [p+1,p+q] : \nabla_i
d(\lambda^{k^*}) \geq 0 \right\}$. Using these notations and the
definition of $\mathbb{D}$, we can write for all $i \in I_{-}$:
\begin{equation}
\label{eq_part2} \left|\left[\nabla_i
d(\lambda^{k^*})\right]_{\rset_+}\right|_{W_{ii}^{-1}}^2=0 \leq
\left| \lambda^{k^*}_i-\lambda^{k^*+1}_i\right|_{W_{ii}}^2.
\end{equation}
On the other hand, for all $i \in I_{+}$ we have:
\begin{align}
\label{eq_part3} \left|\left[\nabla_i
d(\lambda^{k^*})\right]_{\rset_+}\right|_{W_{ii}^{-1}}^2&=\left|\nabla_i
d(\lambda^{k^*})\right|_{W_{ii}^{-1}}^2=\left|\left[W_{ii}^{-1}\nabla_i
d(\lambda^{k^*})\right]_{\rset_+}\right|_{W_{ii}}^2\\
&=\left|\lambda^{k^*}_i-\lambda^{k^*+1}_i\right|_{W_{ii}}^2.
\nonumber
\end{align}
Summing up the relations \eqref{eq_part1},\eqref{eq_part2} and
\eqref{eq_part3} for all $i=1,\dots,p+q$ and combine the result
with \eqref{eq_bound_grad_direction} we obtain:
\begin{equation*}
\left\|\left[\nabla
d(\lambda^{k^*})\right]_{\mathbb{D}}\right\|_{W^{-1}}^2 \leq
\left\|\lambda^{k^*}-\lambda^{k^*+1}\right\|_{W}^2 \leq \frac{4
\|\lambda^*\|_W^2}{(k+1)^3}.
\end{equation*}
Taking now into account that $\left[\nabla
d(\lambda^{k^*})\right]_{\mathbb{D}}=\left[
\begin{array}{c}
  A z^{k^*}-b \\
  ~~\left[C z^{k^*}-c\right]_{\rset^q_+} \\
\end{array}
\right]$ and using the definition of $\mathcal{R}$ we conclude the
result. \qed

\noindent We further characterize  the primal suboptimality and also
the distance from the last iterate $z^{k^*}$ to the optimal solution
$z^*$ of the original optimization problem \eqref{eq_prob_princc}.

\begin{theorem}
\label{theorem_primal_optim_hdwfg} Let the conditions in Theorem
\ref{th_primal_fesability_hdwfg} be satisfied and the function $f$
be Lipschitz continuous with constant $L_{\mathrm{f}}$, i.e.
$|f(z)-f(y)| \leq L_{\mathrm{f}}\|z-y\|$ for all $z,y \in
\rset^n$. Then, the following estimate on primal suboptimality for
problem \eqref{eq_prob_princc} can be derived:
\begin{equation}
\label{bound_primal_optim_hdwfg} -\frac{2
\mathcal{R}^2}{(k+1)\sqrt{(k+1)}}\leq\! f(z^{k^*}) -\! f^*
\!\!\leq
\frac{2L_{\mathrm{f}}\mathcal{R}}{\sqrt{\sigma_{\mathrm{f}}}(k+1)}.
\end{equation}

Moreover, the sequence $z^{k^*}$ converge to the unique optimal
solution $z^*$ of \eqref{eq_prob_princc} with the the following
rate:
\begin{equation}
\label{bound_primal_solution_hdwfg} \|z^{k^*}-z^*\| \leq \frac{2
\mathcal{R}}{\sqrt{\sigma_{\mathrm{f}}}(k+1)}.
\end{equation}
\end{theorem}
\proof The left-hand side inequality of
\eqref{bound_primal_optim_hdwfg} follows using a similar reasoning
as in Theorem \ref{theorem_primal_optim} and the result of Theorem
\ref{th_primal_fesability_hdwfg}. In order to prove the right
hand-side inequality of \eqref{bound_primal_optim_hdwfg} we first
show \eqref{bound_primal_solution_hdwfg}. Using Lemma
\ref{lemma_dist_optim} with $\lambda=\lambda^{k^*}$ we have:
\begin{equation*}
\|z^{k^*}-z^*\| \leq
\sqrt{\frac{2}{\sigma_{\mathrm{f}}}}\sqrt{f^*-d(\lambda^{k^*})}\leq
\frac{2\|\lambda^*\|_W}{\sqrt{\sigma_{\mathrm{f}}}(k+1)},
\end{equation*}
with the last inequality resulting from Theorem
\ref{th_dual_suboptimality_hdfwg}.  Using now the previous
inequality and the Lipschitz property of  $f$ we obtain:
\begin{equation*}
f(z^{k^*})-f^* \leq L_{\mathrm{f}}\|z^{k^*}-z^*\| \leq
\frac{2L_{\mathrm{f}}\|\lambda^*\|_W}{\sqrt{\sigma_{\mathrm{f}}}(k+1)},
\end{equation*}
which together with the definition of $\mathcal{R}$ conclude the
statement. \qed

\begin{remark}
\label{remark_h-dfg} $(i)$ In a similar manner as in Section
\ref{sec_ddfg} using any $\lambda^0 \in \mathbb{D}$ for the
initialization of the algorithm the order of estimates on primal
and dual suboptimality and primal feasibility violation remain the same.\\
$(ii)$ For a given accuracy $\epsilon$, it follows from the results
of this section that we need to perform
$\mathcal{O}\left(\frac{1}{\sqrt[3]{\epsilon^2}}\right)$ iterations
in order to be able to provide a primal suboptimal and
near-feasible solutions based on the last primal iterate $ z^{k^*} $.\\
$(iii)$ Even if the theoretical results show that the estimates on
primal suboptimality and feasibility violation are worse for
Algorithm (\textbf{H-DFG}) in comparison with the ones of Algorithm
(\textbf{DFG}), we have observed that in practice  the last iterate
behaves better. We will discuss this issue in more detail in Section
\ref{sec_numerical}.
\end{remark}

Application of Algorithms (\textbf{DFG}) and (\textbf{H-DFG}) on practical
engineering problems such as (DC-OPF) can be also found in \cite{NecNed:14}.

\section{Linear convergence for dual gradient method under an error bound property}
\label{sec_dg_error_bound}

In this section we show that under the additionally assumption that
the gradients $\nabla f_i$ are Lipchitz continuous the dual problem
\eqref{eq_dual_prob} satisfies a certain error bound property which
allows us to prove a global linear convergence for a distributed
dual gradient method. From our best knowledge this is the first
result showing global  linear convergence of a dual gradient
algorithm. All existing convergence  results from the literature on
dual gradient method   either  show sublinear convergence
\cite{NecSuy:08,NecNed:13,NedOzd:09,BecNed:13} or at most
\textit{local} linear convergence \cite{LuoTse:93a}.

\subsection{Error bound property of the dual problem}

In this section we  assume  that additionally we have Lipschitz
continuity  on the gradient of the primal objective function. Under
strong convexity and this assumption  we prove an error bound type
property on the corresponding dual problem. Our approach for proving
a certain  error bound property is in a way similar to the one in
\cite{LuoTse:92a,WanLin:13}. However, our results are more general
in the sense that we allow the constraints set $\mathbb{D}$ to be an
unbounded polyhedron as opposed to the results in \cite{WanLin:13}
where the authors  show error bound property  only for bounded
polyhedra or the entire space. Also, our  gradient mapping
introduced below is more general than the one used in the standard
analysis of the error bound property (see e.g.
\cite{LuoTse:92a,WanLin:13}). Last but not least important is that
our approach works also for dual problems. Thus, we make further the
following assumption:
\begin{assumption}
\label{ass_lipschitz} The convex functions $f_i$ have Lipschitz
continuous gradients w.r.t. Euclidean norm,  with constants $L_i$.
\end{assumption}
For the convex function $f$, we denote its conjugate by
\cite{RocWet:98}:
\begin{equation*}
\tilde{f}(y)=\sum_{i=1}^M \tilde{f}_i(y),
\end{equation*}
where $\tilde{f}_i(y) =\max \limits_{x_i \in \rset^{n_i}} \langle y,
x_i \rangle -f_i(x_i)$. According to Proposition 12.60 in
\cite{RocWet:98}, under the Assumption \ref{ass_lipschitz} each
function  $\tilde{f}_i(y)$ is strongly convex w.r.t. Euclidean norm,
with constant $\frac{1}{L_i}$, which implies that function
$\tilde{f}$ is strongly convex w.r.t. Euclidean norm, with constant:
\[\sigma_{\tilde{\mathrm{f}}}=\min \limits_{i \in
\left\{1,\dots,M\right\}} \frac{1}L_i.\] Note that in these settings
our dual function can be written as:
\begin{equation}
\label{eq_dual_conj} d(\lambda)=-\tilde{f}(-G^T \lambda)-g^T
\lambda.
\end{equation}

The following lemma whose proof can be also found in \cite[Lemma
4.2]{WanLin:13} or in  \cite[Lemma 3.1]{LuoTse:92a} will help us to
prove the desired error bound property for our dual problem
\eqref{eq_dual_prob}. For completeness we also give the proof.

\begin{lemma}
\label{lemma_unique_tstar} Let Assumptions \ref{ass_strong} and
\ref{ass_lipschitz} hold. Then, there exists a unique $y^* \in
\rset^n$ such that:
\begin{equation}
\label{eq_unique_tstar} G^T \lambda^* =y^* \quad \forall \lambda^*
\in \Lambda^*.
\end{equation}
Moreover, $\nabla d(\lambda) = G \nabla \tilde{f}(-y^*)-g$ is
constant for all $\lambda \in \Lambda$, where the set
$\Lambda=\left\{\lambda \in \mathbb{D} : G^T \lambda =
y^*\right\}$.
\end{lemma}
\proof Let $\lambda_1^*, \lambda_2^* \in \Lambda^*$. From
concavity of $d$ and the fact that the optimal value is the same
for all $\lambda^* \in \Lambda^*$ we have:
\begin{equation*}
d\left(\frac{\lambda_1^*+\lambda_2^*}{2}\right)=
\frac{d(\lambda_1^*)+d(\lambda_2^*)}{2}.
\end{equation*}
Using now \eqref{eq_dual_conj} we can write the following
equality:
\begin{equation*}
-f^*\left(-G^T\frac{\lambda_1^*+\lambda_2^*}{2}\right)
-g^T\frac{\lambda_1^*+\lambda_2^*}{2}=-\frac{f^*(-G^T\lambda_1^*)
+f^*(-G^T\lambda_2^*)}{2}-g^T\frac{\lambda_1^*+\lambda_2^*}{2}.
\end{equation*}
From the strong convexity property of $\tilde{f}$ we have $G^T
\lambda_1^*=G^T \lambda_2^*$, which implies that there exists a
unique $y^*=G^T \lambda^*$ for all $\lambda^* \in \Lambda^*$. The
second statement of the Lemma follows immediately from the
definition of $\Lambda$ and the fact that $y^*$ is unique. \qed

We introduce further the following notations:
\begin{equation}
\label{eq_notations_lambda} r=\left[\lambda\right]_{\Lambda}^W,
\bar{\lambda}=\left[\lambda\right]_{\Lambda^*}^W ~\text{and}
~\bar{r}=\left[r\right]_{\Lambda^*}^W.
\end{equation}

\begin{remark}
\label{remark_nonexpansive}
\begin{itemize}
\item[(i)] Note that for any convex set $D$ and any positive
definite matrix $W$ the projection mapping $[\cdot]_D^W$ is
nonexpansive, i.e. $\|[\lambda]_D^W-[\omega]_D^W\|_W \leq
\|\lambda-\omega\|_W$. In order to prove this nonexpansiveness
property one can use a similar approach as for the Euclidean
projection mapping $[\cdot]_D$. \item[(ii)] In the case $W$ is a
positive definite diagonal matrix and the set $D$ can be written as
the Cartesian product of some sets in $\rset$ we have for any vector
$\lambda \in D$ the following equivalence between projections
$[\lambda]_D^W=[\lambda]_D$. \qed
\end{itemize}
\end{remark}

Using now the notations \eqref{eq_notations_lambda} we can write:
\begin{equation}
\label{eq_intermediate_ineq}  \|\lambda- \bar{\lambda}\|^2_W \leq
\left\|\lambda-\bar{r}\right\|^2_W \leq \left(\|\lambda
-r\|_W+\|r-\bar{r}\|_W\right)^2 \leq
2\|\lambda-r\|_W^2+2\|r-\bar{r}\|^2_W.
\end{equation}
In what follows we will show how we can find upper bounds on
$\|\lambda -r\|_W$ and $\|r-\bar{r}\|_W$ such that we will be able
to establish an error bound property on the dual problem
\eqref{eq_dual_prob}, i.e. there exists a positive constant
$\kappa$, which depends on the original problem  data $L_i$,
$\sigma_i$, $G$ and also on the norm $\|\lambda-\bar{\lambda}\|_W$
such that:
\begin{equation}
\label{eq_ess_err_bound} \|\lambda-\bar{\lambda}\|_W \leq
\kappa(\|\lambda-\bar{\lambda}\|_W) \;  \|\nabla^+ d(\lambda)\|_W
\quad \forall \lambda \in \mathbb{D},
\end{equation}
where the mapping
\begin{equation}
\label{eq_gradient_map} \nabla^+
d(\lambda)=\left[\lambda+W^{-1}\nabla
d(\lambda)\right]_{\mathbb{D}}^W-\lambda
\end{equation}
denotes the  gradient map.   Further, we establish a result which
will help us in proving the error bound \eqref{eq_ess_err_bound}:
\begin{lemma}
\label{lemma_gradient-gradient_map} Let Assumption
\ref{ass_lipschitz} hold and $\nabla^+ d$ be given by
\eqref{eq_gradient_map}. Then, the following inequality holds:
\begin{equation*}
\langle \nabla d(\omega)-\nabla d(\lambda), \lambda -\omega \rangle
\leq 2 \|\nabla^+ d(\lambda)-\nabla^+ d(\omega)\|_W \;
\|\lambda-\omega\|_W \quad \forall \lambda,\omega \in \mathbb{D}.
\end{equation*}
\end{lemma}
\proof First, let us recall that $\left[\lambda+W^{-1}\nabla
d(\lambda)\right]^W_{\mathbb{D}}$ is the unique solution of the
optimization problem:
\begin{equation}
\label{eq_projection_W} \min_{\xi \in \mathbb{D}}
\|\xi-\lambda-W^{-1}\nabla d(\lambda)\|_W^2,
\end{equation}
for which the optimality conditions reads:
\begin{equation*}
\langle W\!\left(\!
\left[\lambda\!+\!W^{-1}\nabla
d(\lambda)\right]^W_{\mathbb{D}}\!\!-\!\left(\!\lambda\!+\!W^{-1}\nabla
d(\lambda)\right)\!\!\right), \xi -\left[\lambda\!+\!W^{-1}\nabla
d(\lambda)\right]^W_{\mathbb{D}} \rangle \geq 0 \quad \forall \xi
\in \mathbb{D}.
\end{equation*}
Taking now $\xi=\left[\omega+W^{-1}\nabla d(\omega
)\right]^W_{\mathbb{D}}$ in the previous  inequality, adding and
substracting $\lambda$ and $\omega$ in the left term of the scalar
product and using the definition of $\nabla^+ d$ we obtain:
\begin{equation*}
\langle W\left(\nabla^+ d(\lambda)-W^{-1}\nabla d(\lambda)\right),
\nabla^+ d(\lambda)+\lambda - \omega -\nabla^+ d(\omega) \rangle
\leq 0,
\end{equation*}
which, together with the fact that $W$ is symmetric implies:
\begin{equation*}
\langle \nabla^+ d(\lambda)-W^{-1}\nabla
d(\lambda),W\left(\lambda-\omega\right)+ W\left(\nabla^+
d(\lambda)-\nabla^+ d(\omega)\right) \rangle \leq 0.
\end{equation*}
Rearranging now the terms in the previous inequality we have:
\begin{align*}
-\langle \nabla d(\lambda),\lambda-\omega\rangle &\leq -\langle
\nabla^+ d(\lambda), W(\lambda-\omega) \rangle + \langle \nabla
d(\lambda), \nabla^+ d(\lambda)-\nabla^+ d(\omega)\rangle\\
&\qquad \qquad - \langle \nabla^+ d(\lambda), W \left(\nabla^+
d(\lambda)- \nabla^+ d(\omega)\right) \rangle.
\end{align*}
Writing now the previous inequality with $\lambda$ and $\xi$
interchanged and summing them up we can write:
\begin{align*}
\langle \nabla &d(\xi)\!-\!\nabla d(\lambda), \lambda\!-\!\xi
\rangle  \!\leq \!\langle
\nabla^+\!d(\xi)\!-\!\nabla^+\!d(\lambda), W(\lambda\!-\!\xi)
\rangle\!\\
&\qquad\qquad\qquad\qquad+\langle \nabla d(\lambda)\!-\!\nabla
d(\xi),
\nabla^+\!d(\lambda)\!-\!\nabla^+\!d(\xi)\rangle\!-\!\|\nabla^+d(\lambda)\!-\!\nabla^+
d(\xi)\|_W^2\\
&\qquad\leq \langle \nabla^+\!d(\xi)\!-\!\nabla^+\!d(\lambda),
W(\lambda\!-\!\xi) \rangle\! +\langle \nabla d(\lambda)\!-\!\nabla
d(\xi), \nabla^+\!d(\lambda)-\nabla^+ d(\xi) \rangle\\
&\qquad\leq \|\nabla^+ d(\lambda)-\nabla^+ d(\xi)
\|_W\left(\|W(\lambda-\xi)\|_{W^{-1}}+\|\nabla d(\lambda)-\nabla
d(\xi)\|_{W^{-1}}\right)\\
&\qquad\leq 2\|\nabla^+ d(\lambda)-\nabla^+ d(\xi) \|_W
\|\lambda-\xi\|_W,
\end{align*}
which concludes the statement. \qed

The next lemma gives un upper bound on $\|\lambda -r\|_W$:
\begin{lemma}
\label{lemma_upper_r} Under the Assumptions \ref{ass_strong}  and
\ref{ass_lipschitz} there exists a positive constant $\kappa_1$ such
that the following inequality holds:
\begin{equation}
\label{eq_bound_r} \|\lambda-r\|^2_W \leq
 \kappa_1 \; \|\nabla^+d(\lambda)\|_W  \; \|\lambda-\bar{\lambda}\|_W \quad \forall \lambda
\in \mathbb{D},
\end{equation}
where $\kappa_1=\frac{2}{\sigma_{\tilde{\mathrm{f}}}}\theta_1^2$,
with $\theta_1$ being a finite constant depending on the matrix
$G$.
\end{lemma}
\proof  First, let us notice that we can write the set $\Lambda$
as the following set characterized by linear equalities and
inequalities:
\begin{align}
\label{eq_sets_mangasarian2} \Lambda&=\left\{\omega \in
\rset^{p+q}:~F \omega \leq 0,~~ G^T \omega = y^* \right\},
\end{align}
where $F=\left[0_{q,p} -I_q\right]$. Since $\lambda \in
\mathbb{D}$ this implies $F \lambda \leq 0$ and therefore,
according to Theorem 2 in \cite{Rob:73} we can write:
\begin{equation}
\label{eq_mangasarian} \|\lambda-r\|_W \leq \theta_1 \|G^T \lambda
- y^*\|_{\infty} \leq \theta_1 \|G^T \lambda - y^*\|,
\end{equation}
where $\theta_1$ is finite and depends only on the matrix $G$ and
on the norms $\|\cdot\|_W$ and $\|\cdot\|_{\infty}$.
From the strong convexity property of $\tilde{f}$ combined with
the fact that $G^T \bar{\lambda}=y^*$ we have:
\begin{align}
\label{eq_strong_fstar_bound} \sigma_{\tilde{\mathrm{f}}}\|G^T
\lambda -y^*\|^2 &\leq \langle \nabla \tilde{f}(-G^T
\lambda)-\nabla \tilde{f}(-G^T
\bar{\lambda}), -G^T \lambda +G^T \bar{\lambda} \rangle \nonumber\\
&=\langle -G\nabla \tilde{f}(-G^T \lambda)+g +G\nabla
\tilde{f}(-G^T
\bar{\lambda})-g, \lambda - \bar{\lambda} \rangle  \\
&=\langle \nabla d(\bar{\lambda}) -\nabla d(\lambda),
\lambda-\bar{\lambda}\rangle \nonumber \\
&\leq 2\|\nabla^+ d(\lambda) - \nabla^+
d(\bar{\lambda})\|_W\|\lambda-\bar{\lambda}\|_W=2\|\nabla^+
d(\lambda)\|_W\|\lambda-\bar{\lambda}\|_W, \nonumber
\end{align}
where the last inequality follows from Lemma
\ref{lemma_gradient-gradient_map} with $\lambda=\lambda$ and
$\xi=\bar{\lambda}$ and the last equality is deduced from the fact
that since $\bar{\lambda} \in \Lambda^*$ this implies that
$\nabla^+ d(\bar{\lambda})=0$. Combining now
\eqref{eq_mangasarian} with \eqref{eq_strong_fstar_bound} we
obtain the result. \qed

The following result establishes an upper bound on
$\|r-\bar{r}\|_W$:

\begin{lemma}
\label{lemma_upper_ru} Let Assumptions \ref{ass_strong} and
\ref{ass_lipschitz} hold. Then, the following inequality is valid:
\begin{equation}
\label{eq_upper_ru} \|r-\bar{r}\|_W^2 \leq \kappa_2 \|\nabla^+
d(\lambda)\|_W\|\lambda-\bar{\lambda}\|_W \quad \forall \lambda \in
\mathbb{D},
\end{equation}
where $r, \bar{r}$ are given by \eqref{eq_notations_lambda} and
\[\kappa_2=6\theta_2^2\left(2\max \limits_{\lambda^* \in
\Lambda^*}\|\lambda-\lambda^*\|_W^2+2\|\nabla
d(\bar{\lambda})\|_{W^{-1}}^2\right)\left(1+3\theta_1^2\frac{2}{\sigma_{\tilde{\mathrm{f}}}}\right),\]
with $\theta_2$ is a constant depending on $C$, $\nabla
d(\bar{\lambda})$ and $y^*$.
\end{lemma}
\proof Since $\Lambda^* \subseteq \Lambda \subseteq \mathbb{D}$
and $G^T \xi = y^*$ for all $\xi \in \Lambda$, then the dual
problem \eqref{eq_dual_prob} has the same optimal solutions as the
following linear problem:
\begin{equation}
\label{eq_linear1} \arg \max_{\lambda \in \mathbb{D}}
d(\lambda)=\arg \max_{\xi \in \Lambda} d(\xi)= \arg \max_{\xi \in
\Lambda} -\tilde{f}(-y^*) -\langle g,\xi\rangle =\arg \max_{\xi
\in \Lambda} -\langle g,\xi\rangle.
\end{equation}
Further, let us recall that $\nabla d(\bar{\lambda})= G
\tilde{f}(-y^*)-g$ for any $\lambda \in \mathbb{D}$ and thus we
have that $\langle \nabla d(\bar{\lambda}), \xi\rangle= -\langle
\nabla \tilde{f}(-y^*), y^*\rangle-\langle g, \xi\rangle$ for all
$\xi \in \Lambda$. Therefore, we can write further:
\begin{equation}
\label{eq_linear2} \arg \max_{\xi \in \Lambda} \langle \nabla
d(\bar{\lambda}), \xi\rangle =\arg \max_{\xi \in \Lambda}- \langle
\nabla \tilde{f}(-y^*), y^*\rangle -\langle g, \xi\rangle=\arg
\max_{\xi \in \Lambda} -\langle g, \xi\rangle.
\end{equation}
Combining now \eqref{eq_linear1} with \eqref{eq_linear2} we can
conclude that any solution $\bar{\xi}=[\xi]_{\Lambda^*}^W$ with
$\xi \in \Lambda$ of the dual problem \eqref{eq_dual_prob} is also
a solution of problem \eqref{eq_linear2}. Since $\nabla_{\nu}
d(\bar{\lambda})=Az^*-b = 0$ and $\nabla_{\mu}
d(\bar{\lambda})=Cz^*-c \leq 0$, then we also have that the
maximum in \eqref{eq_linear2} is finite and thus problem
\eqref{eq_linear2} is solvable. Applying now Theorem 2 in
\cite{Rob:73} to the optimality conditions of problem
\eqref{eq_linear2} and its dual we obtain:
\begin{equation}
\label{eq_robinson} \|\xi-\bar{\xi}\|_W \leq \theta_2 |\langle
\nabla d(\bar{\lambda}), \xi \rangle -\langle \nabla
d(\bar{\lambda}), \bar{\xi}\rangle| \quad \forall \xi \in \Lambda,
\end{equation}
where $\theta_2$ is a constant depending only on the matrix $C$
and on vectors $\nabla d(\bar{\lambda})$ and $y^*$ (see eq. (6) in
\cite{Rob:73} for details). Using the previous relation we have:
\begin{equation}
\label{eq_linear3} \|\xi-\bar{\xi}\|_W \leq \theta_2 |\langle
\nabla d(\bar{\lambda}), \xi \rangle -\langle \nabla
d(\bar{\lambda}), \bar{\xi}\rangle|=\theta_2 \langle \nabla
d(\bar{\lambda}), \bar{\xi}-\xi \rangle.
\end{equation}
For any $\xi \in \Lambda$ the optimality conditions of the
following  projection problem $\min \limits_{\omega \in
\Lambda}\|\omega - \xi - W^{-1}\nabla d(\bar{\lambda}) \|_W^2$
become:
\begin{equation*}
\langle W\left(\left[\xi+W^{-1}\nabla
d(\bar{\lambda})\right]_{\Lambda}^W-\xi-W^{-1}\nabla
d(\bar{\lambda})\right),\left[\xi+W^{-1}\nabla
d(\bar{\lambda})\right]_{\Lambda}^W- \omega \rangle \leq 0
\end{equation*}
for all $\omega \in \Lambda$. Taking now $\omega=
\bar{\xi}=[\xi]^W_{\Lambda^*}$ and since $W$ is a symmetric matrix
we obtain:
\begin{align*}
\langle\nabla
d&(\bar{\lambda}),\bar{\xi}\!-\!\xi\rangle\leq\langle
\left[\xi\!+\!W^{-1}\nabla
d(\bar{\lambda})\right]_{\Lambda}^W\!\!-\xi, W\left(
\bar{\xi}-\left[\xi+W^{-1}\nabla
d(\bar{\lambda})\right]_{\Lambda}^W\right)+\nabla d(\bar{\lambda}) \rangle\\
&= \!\langle \left[\xi\!+\!W^{-1}\nabla
d(\bar{\lambda})\right]_{\Lambda}^W\!\!-\xi,W\!\left(
\xi\!-\!\left[\xi\!+\!W^{-1}\nabla
d(\bar{\lambda})\right]_{\Lambda}^W\right)+W(\bar{\xi}\!-\!\xi)\!+\!\nabla d(\bar{\lambda})\rangle\\
&\leq \langle \left[\xi+W^{-1}\nabla
d(\bar{\lambda})\right]_{\Lambda}^W-\xi,W(\bar{\xi}-\xi)+\nabla
d(\bar{\lambda})\rangle\\
&\leq \left\|\left[\xi+W^{-1}\nabla
d(\bar{\lambda})\right]_{\Lambda}^W-\xi\right\|_W\left(\|W(\bar{\xi}-\xi)\|_{W^{-1}}+\|\nabla
d(\bar{\lambda})\|_{W^{-1}}\right)\\
&=\|\nabla^{+} d(\xi)\|_W\left(\|\xi-\bar{\xi}\|_{W}+\|\nabla
d(\bar{\lambda})\|_{W^{-1}}\right),
\end{align*}
where in the last equality we used the definition of $\nabla^+ d$
and the fact that $\nabla d(\bar{\lambda})= \nabla d(\xi)$ for all
$\xi \in \Lambda$ (see Lemma \ref{lemma_unique_tstar}). Combining
now the previous inequality with \eqref{eq_linear3}, taking $\xi=r
\in \Lambda$ and squaring both sides we obtain:
\begin{equation}
\label{eq_bound_r_inter} \|r-\bar{r}\|_W^2 \leq \theta_2^2
\left(\|r-\bar{r}\|_W+\|\nabla
d(\bar{\lambda})\|_{W^{-1}}\right)^2\|\nabla^+ d(r)\|_W^2.
\end{equation}
Since $\bar{r}=\left[r\right]_{\Lambda^*}^W$ and $\Lambda^*
\subseteq \Lambda$ we also have
$\bar{r}=\left[\bar{r}\right]_{\Lambda}^W$. Thus, using the
nonexpansive  property of the projection we can write:
\begin{equation}
\label{eq_bound_r_max} \|r-\bar{r}\|_W \leq \|\lambda- \bar{r}\|_W
\leq \max \limits_{\lambda^* \in
\Lambda^*}\|\lambda-\lambda^*\|_W.
\end{equation}
Further, our goal is to find an upper bound for $\|\nabla^+
d(r)\|_W$ in terms of $\|\lambda-\bar{\lambda}\|_W$ and
$\|\nabla^+ d(\lambda)\|_W$. For this purpose let us first prove
that $\nabla^+ d$ is Lipschitz continuous with constant $3$ w.r.t
to the norm $\|\cdot\|_W$. For any $\lambda, \tilde{\lambda} \in
\mathbb{D}$ we can write:
\begin{align}
\label{eq_mapping_lipsch} \|\nabla^+ d(\lambda)\!-\!\nabla^+
d(\tilde{\lambda})\|_W &\leq \|\lambda\!-\!
\tilde{\lambda}\|_W+\left\|\left[\lambda\!+\!W^{-1}\nabla
d(\lambda)\right]_{\mathbb{D}}^W-\left[\tilde{\lambda}\!+\!W^{-1}\nabla
d(\tilde{\lambda})\right]_{\mathbb{D}}^W\right\| \nonumber\\
&\leq \|\lambda-\tilde{\lambda}\|_W+\|\lambda+W^{-1}\nabla
d(\lambda)-\tilde{\lambda}-W^{-1}\nabla d(\tilde{\lambda})\|_W \\
&\leq 2\|\lambda - \tilde{\lambda}\|_W+\|\nabla d(\lambda)-\nabla
d(\tilde{\lambda})\|_{W^{-1}} \leq 3\|\lambda-\tilde{\lambda}\|_W.
\nonumber
\end{align}
Using now \eqref{eq_mapping_lipsch} with
$\tilde{\lambda}=\bar{\lambda}$ and taking into account that
$\nabla^+ d(\bar{\lambda})=0$, we have:
\begin{equation}
\label{eq_lemma_a5} \|\nabla^+ d(\lambda)\|_W=\|\nabla^+
d(\lambda)-\nabla^+ d(\bar{\lambda})\|_W \leq
3\|\lambda-\bar{\lambda}\|_W.
\end{equation}
Using now again \eqref{eq_mapping_lipsch} and the Lipschitz
continuity property of $\nabla d$ we can write:
\begin{align}
\label{eq_ineq_error_bound}  \|\nabla^+ d(r)\|_W^2 &\leq
\left(\|\nabla^+ d(\lambda)\|_W+\|\nabla^+
d(r)-\nabla^+ d(\lambda)\|_W\right)^2 \nonumber \\
& \leq 2\|\nabla^+ d(\lambda)\|_W^2+2\|\nabla^+ d(r)-\nabla^+
d(\lambda)\|_W^2 \nonumber\\
&\leq 6\|\nabla^+
d(\lambda)\|_W\|\lambda-\bar{\lambda}\|_W+18\|\lambda-r\|_W^2 \nonumber\\
&\leq
6\left[1+3\theta_1^2\frac{2}{\sigma_{\tilde{\mathrm{f}}}}\right]\|\nabla^+
d(\lambda)\|_W\|\lambda-\bar{\lambda}\|_W,
\end{align}
where in the last inequality we used \eqref{eq_bound_r}. Introducing
now \eqref{eq_bound_r_max} and \eqref{eq_ineq_error_bound} in
\eqref{eq_bound_r_inter} and using the inequality $(\alpha+\beta)^2
\leq 2\alpha^2+2\beta^2$ we obtain the result. \qed

Note that for any finite $\lambda \in \mathbb{D}$, since $\Lambda^*$
is bounded we have that $\max \limits_{\lambda^* \in \Lambda^*}
\|\lambda-\lambda^*\|_W$ is also finite. The following theorem
establishes an error bound type property for the dual problem
\eqref{eq_dual_prob}:

\begin{theorem}
\label{theorem_error_bound} Let Assumptions \ref{ass_strong}  and
\ref{ass_lipschitz} hold. Then, there exists $\kappa$, depending on
the data of the original problem and $\max \limits_{\lambda^* \in
\Lambda^*}\|\lambda - \lambda^*\|_W$, such that the following error
bound property can be established for the dual problem
\eqref{eq_dual_prob}:
\begin{equation}
\label{eq_error_bound} \|\lambda - \bar{\lambda}\|_W \leq
\kappa\left(\|\lambda - \lambda^*\|_W \right) \|\nabla^+
d(\lambda)\|_W ~~ \forall \lambda \in \mathbb{D},
\end{equation}
with:
\begin{align*}
\kappa&\left(\|\lambda\!-\!\lambda^*\|_W
\right)\!=\!\kappa_1\!+\!\kappa_2\!=\!\theta_1^2
\frac{4}{\sigma_{\tilde{\mathrm{f}}}}+12\theta_2^2\!\left(\!\!2\!\max
\limits_{\lambda^* \in
\Lambda^*}\|\lambda\!-\!\lambda^*\|_W^2\!+\!2\|\nabla
d(\bar{\lambda})\|_W^2\!\!\right)\!\!\left(\!1\!+\!3\theta_1^2\frac{2}{\sigma_{\tilde{\mathrm{f}}}}\right).
\end{align*}
\end{theorem}
\proof The result follows immediately by using \eqref{eq_bound_r}
from Lemma \ref{lemma_upper_r} and \eqref{eq_upper_ru} from Lemma
\ref{lemma_upper_ru} in \eqref{eq_intermediate_ineq} and dividing
both sides by $\|\lambda -\bar{\lambda}\|_W$. \qed

\subsection{Convergence analysis using the last iterate}
\label{subsec_convergence_dfg_error_bound} Under the error bound
property for the dual problem defined in Theorem
\ref{theorem_error_bound}  we will show in this section   linear
convergence for a  distributed dual gradient method.    From our
knowledge this  is the  first result showing  \textit{global} linear
convergence rate on primal suboptimality and infeasibility for the
last primal iterate of a dual gradient algorithm, as opposed  to the
results in \cite{LuoTse:93a} where only \textit{local} linear
convergence was derived for such an algorithm. Thus, we now
introduce the following distributed dual gradient method:
\begin{center}
\framebox{
\parbox{8.3cm}{
\begin{center}
\textbf{ Algorithm {\bf (DG)}}
\end{center}
{Initialization: $\lambda^0=0$. For $k\geq 0$ compute:}
\begin{enumerate}
\item{$z^k = \arg \min\limits_{z \in \rset^n}
\mathcal{L}(z,\lambda^k)$}.

\item ${\lambda}^{k+1} = \left[\lambda^k+W^{-1}\nabla
d(\lambda^k)\right]_{\mathbb{D}}$.
\end{enumerate}
}}
\end{center}

The next lemma, which is a generalization of  a known result for the
gradient method (see e.g. \cite{Nes:04,NecNed:13}) will help us to
analyze the convergence of the Algorithm (\textbf{DG}):

\begin{lemma}
\label{lemma_decrese_dist_dg} Let Assumption \ref{ass_strong}
hold and the sequence $\{\lambda^k\}_{k \geq 0}$ be generated by
Algorithm (\textbf{DG}). Then, the following inequalities hold:
\begin{equation}
\label{eq_decrease} \|\lambda^k-\lambda^*\|_W \leq \cdots \leq
\|\lambda^0-\lambda^*\|_W \quad \forall \lambda^* \in \Lambda^*, k
\geq 0.
\end{equation}
\end{lemma}
\proof Taking $\lambda=\lambda^*$ in the optimality condition of
\eqref{eq_projection_dual}, we obtain the following inequality:
\begin{equation}
\label{eq_opt_cond_decrease} \langle \nabla d(\lambda^k)-W
(\lambda^{k+1}-\lambda^k),\lambda^*-\lambda^{k+1}\rangle \leq 0.
\end{equation}
Further, we can write:
\begin{align}
\label{inequalities_gradient}
\|\lambda^{k+1}\!\!-\!\lambda^*\|_W^2\!
&=\!\|\lambda^{k+1}-\lambda^k+\lambda^k-\lambda^*\|_W^2 \nonumber\\
&=\!\|\lambda^k\!\!-\!\lambda^*\|_W^2\!+\!2\langle W(
\lambda^{k+1}\!\!-\!\lambda^k),\lambda^{k}\!\!-\!\lambda^{k+1}\!\!+\!\lambda^{k+1}\!\!-\!\lambda^*
\rangle\!+\!\|\lambda^{k+1}\!\!-\!\lambda^k\|_W^2 \nonumber\\
&=\!\|\lambda^k-\lambda^*\|^2_W+2\langle W(
\lambda^{k+1}-\lambda^k),\lambda^{k+1}-\lambda^*
\rangle-\|\lambda^{k+1}-\lambda^k\|_W^2 \nonumber \\
&\!\leq \|\lambda^k-\lambda^*\|^2_W-2\langle \nabla
d(\lambda^k),\lambda^*-\lambda^k\rangle \\
&\qquad+2\left(\langle \nabla
d(\lambda^k),\lambda^{k+1}-\lambda^k\rangle\!-\frac{1}{2}
\|\lambda^{k+1}\!\!-\lambda^k\|^2_W\right)  \nonumber\\
&\!\leq
\|\lambda^k-\lambda^*\|^2_W+\!2\left(d(\lambda^k)\!-\!d(\lambda^*)\right)\!+\!2\left(d(\lambda^{k+1})-d(\lambda^k)\right)
\nonumber \\
&\!=\|\lambda^k-\lambda^*\|^2_W+2\left(d(\lambda^{k+1})-d(\lambda^*)\right)
\leq \|\lambda^k-\lambda^*\|^2_W,
\end{align}
where the first inequality follows from
\eqref{eq_opt_cond_decrease} and the second one is derived from
the concavity of the function $d$ and Lemma \ref{descent_lemma}.
Applying now recursively the previous inequality we obtain
\eqref{eq_decrease}. \qed

Using now \eqref{eq_decrease} and the definition of $\mathcal{R}$ we
can write for all $k \geq 0$:
\begin{equation}
\label{general_mult_bound} \max \limits_{\lambda^* \in
\Lambda^*}\|\lambda^{k}-\lambda^*\|_W \leq \max \limits_{\lambda^*
\in \Lambda^*}\|\lambda^0-\lambda^*\|_W=\mathcal{R},
\end{equation}
where the last equality follows from the definition of
$\mathcal{R}$ and the fact that $\lambda^0=0$. Introducing now
this inequality in \eqref{eq_error_bound} we obtain:
\begin{equation}
\label{global_error_bound} \|\lambda^k-\bar{\lambda}^k\|_W \leq
\overline{\kappa} \;  \|\nabla^+ d(\lambda^k)\|_W \quad \forall k
\geq 0,
\end{equation}
where:
\begin{align*}
\overline{\kappa}= \theta_1^2
\frac{4}{\sigma_{\tilde{\mathrm{f}}}}+12\theta_2^2\left(2\mathcal{R}^2+2\|\nabla
d(\bar{\lambda}^k)\|_W^2\right)\left(1+3\theta_1^2
\frac{2}{\sigma_{\tilde{\mathrm{f}}}}\right).
\end{align*}

Using now Remark \ref{remark_nonexpansive} $(ii)$ and the
definition of $\nabla^+ d$ we have:
\begin{equation}
\label{eq_decr_map} \|\nabla^+ d(\lambda^k)\|_W
=\|\lambda^{k+1}-\lambda^k\|_W.
\end{equation}
Further, combining \eqref{global_error_bound} with
\eqref{eq_decr_map} we can write:
\begin{equation}
\label{eq_relation1} \|\lambda^k-\bar{\lambda}^k\|_W \leq
\overline{\kappa}\|\nabla^+
d(\lambda^k)\|_W=\overline{\kappa}\|\lambda^{k+1}-\lambda^k\|_W.
\end{equation}
The following theorem provides an estimate on the dual
suboptimality for Algorithm ({\bf DG}):
\begin{theorem}
\label{theorem_dual_optim_dg} Let Assumptions \ref{ass_strong} and
\ref{ass_lipschitz} hold and the sequences
$\left(z^k,\lambda^k\right)_{k\geq 0}$ be generated by algorithm
({\bf DG}). Then, an estimate on dual suboptimality for
\eqref{eq_dual_prob} is given by:
\begin{equation}
\label{bound_dual_optim_dg} f^*-d({\lambda}^{k+1})\leq
\left(\frac{4(1+\overline{\kappa})}{1+4(1+\overline{\kappa})}\right)^k
\left( f^*-d(\lambda^0)\right).
\end{equation}
\end{theorem}
\proof

First, let us notice that for any $k  \geq 0$, $\lambda^{k+1}$ can
be computed as the unique optimal solution of problem
\eqref{eq_projection_dual}. Thus, from the optimality condition of
problem \eqref{eq_projection_dual} we have:
\begin{equation}
\label{eq_optimality_lambdak} \langle \nabla
d(\lambda^k),\bar{\lambda}^k-\lambda^{k+1}\rangle \leq \langle W
(\lambda^{k+1}-\lambda^k), \bar{\lambda}^k-\lambda^{k+1}\rangle
\leq 0,
\end{equation}
where we recall that $\bar{\lambda}^k=[\lambda^k]^W_{\Lambda^*}$.
Further, since the optimal value of the dual function is unique we
can write:
\begin{align}
\label{ineq_linear_dg}
f^*-d(\lambda^{k+1})&=d(\bar{\lambda}^k)-d(\lambda^{k+1})\leq
\langle \nabla d(\lambda^{k+1}), \bar{\lambda}^k - \lambda^{k+1}
\rangle \nonumber\\
& =\langle \nabla d(\lambda^{k+1})-\nabla d(\lambda^k),
\bar{\lambda}^k - \lambda^{k+1} \rangle +\langle \nabla
d(\lambda^{k}), \bar{\lambda}^k - \lambda^{k+1} \rangle\\
& \leq \| \nabla d(\lambda^{k+1})-\nabla d(\lambda^k)\|_{W^{-1}}
\|\bar{\lambda}^k - \lambda^{k+1}\|_{W} \nonumber\\
&\qquad \qquad \qquad+ \langle W (\lambda^{k+1}-\lambda^k),
\bar{\lambda}^k-\lambda^{k+1}\rangle \nonumber\\
& \leq \|\lambda^{k+1}-\lambda^k\|_{W} \|\bar{\lambda}^k -
\lambda^{k+1}\|_{W} \nonumber+ \|\lambda^{k+1}-\lambda^k\|_{W}
\|\bar{\lambda}^k-\lambda^{k+1}\|_{W} \nonumber\\
&=2\|\lambda^{k+1}-\lambda^k\|_{W}\|\bar{\lambda}^k-\lambda^{k+1}\|_{W}
\nonumber,
\end{align}
where the second inequality follows from
\eqref{eq_optimality_lambdak}. Using now relation
\eqref{eq_relation1} we can write:
\begin{align*}
\|\bar{\lambda}^k-\lambda^{k+1}\|_{W} \leq
\|\bar{\lambda}^k-\lambda^{k}\|_{W}+\|{\lambda}^k-\lambda^{k+1}\|_{W}\leq
\left(1+\overline{\kappa}\right)\|{\lambda}^k-\lambda^{k+1}\|_{W}.
\end{align*}
Introducing now the previous inequality in \eqref{ineq_linear_dg}
and using Lemma \ref{lemma_dual_ascent} we have:
\begin{align*}
f^*-d(\lambda^{k+1})\leq
2\left(1+\overline{\kappa}\right)\|{\lambda}^k-\lambda^{k+1}\|_{W}^2
\leq 4\left(1+\overline{\kappa}\right)
\left(d(\lambda^{k+1})-d(\lambda^k)\right).
\end{align*}
Rearranging the terms in the previous inequality we obtain:
\begin{equation}
\label{eq_recursion} f^*-d(\lambda^{k+1}) \leq
\frac{4(1+\overline{\kappa})}{1+4(1+\overline{\kappa})}\left(f^*-d(\lambda^k)\right).
\end{equation}
Applying now \eqref{eq_recursion} recursively we obtain
\eqref{bound_dual_optim_dg}. \qed

In order to characterize the dual suboptimality we extend the proof for the centralized
gradient algorithm \cite{WanLin:13,LuoTse:93a} to the case of distributed dual gradient Algorithm (\textbf{DG}).
The following theorem gives an estimate on the primal feasibility
violation for Algorithm (\textbf{DG}):

\begin{theorem}
\label{theorem_primal_fesa_dg} Under the assumptions of Theorem
\ref{theorem_dual_optim_dg}, the following estimate holds for the
primal feasibility violation:
\begin{equation}
\label{eq_primal_feasibility_dg} \left\|\left[
\begin{array}{c}
  A z^k-b \\
  ~~\left[C z^k-c\right]_{\rset^q_+} \\
\end{array}%
\right]\right\|_{W^{-1}} \leq
\left(\frac{4(1+\overline{\kappa})}{1+4(1+\overline{\kappa})}\right)^{\frac{k-1}{2}}
\sqrt{2\left( f^*-d(\lambda^0)\right)}.
\end{equation}
\end{theorem}
\proof Using the descent property of dual gradient method
\eqref{eq_ascent_dual} we have:
\begin{align}
\label{eq_feasibility1_dg} \|\lambda^{k}-\lambda^{k+1}\|_W^2 &\leq
2 \left(d(\lambda^{k+1})-d(\lambda^k)\right) \leq
2\left(f^*-d(\lambda^k)\right)\\
&\leq
2\left(\frac{4(1+\overline{\kappa})}{1+4(1+\overline{\kappa})}\right)^{k-1}
\left( f^*-d(\lambda^0)\right), \nonumber
\end{align}
where in the last inequality we used Theorem
\ref{theorem_dual_optim_dg}. Using now a similar reasoning as in
Theorem \ref{th_primal_fesability_hdwfg}, we obtain:
\begin{equation*}
\left\|\left[\nabla
d(\lambda^{k})\right]_{\mathbb{D}}\right\|_{W^{-1}}^2 \leq
\|\lambda^{k}-\lambda^{k+1}\|_W^2 \leq
2\left(\frac{4(1+\overline{\kappa})}{1+4(1+\overline{\kappa})}\right)^{k-1}
\left( f^*-d(\lambda^0)\right),
\end{equation*}
where in second inequality we used \eqref{eq_feasibility1_dg}.
Squaring now both sides of previous inequality and taking into
account the definitions of $\nabla  d$ and $\mathbb{D}$ we obtain
\eqref{eq_primal_feasibility_dg}. \qed

We now characterize  the primal suboptimality and the distance
from the last iterate $z^{k}$ generated by Algorithm (\textbf{DG})
to the optimal solution $z^*$ of our original optimization problem
\eqref{eq_prob_princc}.

\begin{theorem}
\label{theorem_primal_optim_dg} Let the conditions in Theorem
\ref{theorem_primal_fesa_dg} be satisfied. Then, the following
estimate on primal suboptimality for problem
\eqref{eq_prob_princc} can be derived:
\begin{equation}
\label{bound_primal_optim_dg}
-\left(\frac{4(1+\overline{\kappa})}{1+4(1+\overline{\kappa})}\right)^{\frac{k-1}{2}}\mathcal{R}
\sqrt{2\left( f^*-d(\lambda^0)\right)} \leq f(z^k)-f^* \leq v(k),
\end{equation}
where
\begin{align*}
v(k)&=\frac{\mathcal{R}}{\underline{w}}\|G\|\left(\frac{4(1+\overline{\kappa})}{1+4(1+\overline{\kappa})}\right)^{\frac{k}{2}}
\sqrt{\frac{2}{\sigma_{\text{f}}}\left(
f^*-d(\lambda^0)\right)}\\
&\qquad+\frac{\max\limits_{i=1,\dots,M}L_i}{2}
\left(\frac{4(1+\overline{\kappa})}{1+4(1+\overline{\kappa})}\right)^{k}
\frac{2}{\sigma_{\text{f}}}\left( f^*-d(\lambda^0)\right).
\end{align*}
Moreover, the sequence $z^{k}$ converge to the unique optimal
solution $z^*$ of \eqref{eq_prob_princc} with the the following
rate:
\begin{equation}
\label{bound_primal_solution_dg} \|z^{k}-z^*\| \leq
\left(\frac{4(1+\overline{\kappa})}{1+4(1+\overline{\kappa})}\right)^{\frac{k}{2}}
\sqrt{\frac{2}{\sigma_{\text{f}}}\left( f^*-d(\lambda^0)\right)}.
\end{equation}
\end{theorem}
\proof The left-hand side inequality of
\eqref{bound_primal_optim_dg} follows using a similar reasoning as
in Theorem \ref{theorem_primal_optim} and the result of Theorem
\ref{theorem_primal_fesa_dg}. In order to prove the right hand-side
inequality of \eqref{bound_primal_optim_dg} we first show
\eqref{bound_primal_solution_dg}. Using Lemma \ref{lemma_dist_optim}
with $\lambda=\lambda^{k}$ we have:
\begin{equation*}
\|z^{k}-z^*\| \leq
\sqrt{\frac{2}{\sigma_{\mathrm{f}}}}\sqrt{f^*-d(\lambda^{k^*})}\leq
\left(\frac{4(1+\overline{\kappa})}{1+4(1+\overline{\kappa})}\right)^{\frac{k}{2}}
\sqrt{\frac{2}{\sigma_{\text{f}}}\left( f^*-d(\lambda^0)\right)},
\end{equation*}
with the last inequality resulting from Theorem
\ref{theorem_dual_optim_dg}. Let us introduce further the notation
by $\underline{w}=\lambda_{\text{min}}(W)$. Using now the
continuous Lipschitz property of $\nabla f$ we obtain:
\begin{align*}
f(z^{k})-f^* &\leq \langle \nabla f(z^*), z^k-z^*
\rangle+\frac{\max_i L_i}{2}\|z^{k}-z^*\|^2\\
&=\langle -G^T \lambda^*, z^k-z^* \rangle+\frac{\max_i L_i}{2}\|z^{k}-z^*\|^2\\
&\leq
\frac{\mathcal{R}\|G\|}{\underline{w}}\|z^{k}-z^*\|+\frac{\max_i
L_i}{2}\|z^{k}-z^*\|^2,
\end{align*}
where the first equality is deduced from the optimality conditions
of problem $z^*=\arg\min f(z)+\langle \lambda^*, Gz-g\rangle$ and
in the last inequality we used Cauchy-Schwartz inequality, the
fact that $\|\cdot\| \leq \frac{1}{\underline{w}}\|\cdot\|_W$ and
the definition of $\mathcal{R}$. Using now
\eqref{bound_primal_solution_dg} in the previous inequality we
obtain the result. \qed

\section{Distributed implementation}
\label{subsec_distributed_dfwg}

In this section we analyze the distributed implementation of
Algorithms (\textbf{DFG}), (\textbf{H-DFG}) and (\textbf{DG}). We
look first at step $1$ of the Algorithm (\textbf{DFG}). Note that
this step is similar with the steps $1$ of phases $1$ and $2$ of
Algorithm (\textbf{H-DFG}) and the step $1$ of Algorithm
(\textbf{DG}) and therefore their analysis follows in a similar
way. According to \eqref{eq_separable_dual}, for all $i \in V_1$
we have:

\begin{align}
\label{distributed_primal_update} z_i^{k}&=\arg\min_{z_i \in
\rset^{n_i}} f_i(z_i) +\left\langle \lambda^k, \left[A_{ i}^T
C_{ i}^T\right]^Tz_i \right\rangle \nonumber\\
&=\arg \min_{z_i \in \rset^{n_i}} f_i(z_i) +\sum_{j \in \Ncal_i}
\left(\left[A_{ji}^T C_{ji}^T\right]\lambda_j^k\right)^T z_i.
\end{align}
Thus, in order to compute $z_i^k$ the algorithm requires only
local information, namely
$\left\{A_{ji},C_{ji},\lambda^k_j\right\}_{j \in \Ncal_i}$. For
example, in the case of (NUM) problem, the update of source rate
$z_i^k$ requires only the link prices $\lambda^k_j$ which are
utilized by source $i$. Using now the definitions of $W$ and
$\nabla d$, step $2$ in Algorithm (\textbf{DFG}) can be written in
the following form:
\begin{equation}
\label{distributed_dual_update}
\hat{\lambda}_j^k=\left[\lambda_j^k+\left[\begin{array}{c}
   W_{\nu jj}^{-1} \sum_{i \in \bar{\Ncal}_j} A_{ji} z_i^k\\
   W_{\mu jj}^{-1} \sum_{i \in \bar{\Ncal}_j} C_{ji} z_i^k\\
\end{array}\right]\right]_{\rset^{p_j}\times \rset^{q_j}_+},
\forall j \in V_2,
\end{equation}
where $W_{\nu jj}$ and $W_{\mu jj}$ denote the $j$th
block-diagonal element of matrix $W_{\nu}$ and $W_{\mu}$,
respectively. Taking into account the definitions of $W_{\nu jj}$
and $W_{\mu jj}$ we can conclude that in order to update the dual
variable $\hat{\lambda}_j^k$ in step $2$ of Algorithm
(\textbf{DFG}) we require only local information
$\left\{L_{\text{d}_i},A_{ji},C_{ji},z_i^k\right\}_{i \in
\bar{\Ncal}_j}$. Thus, in the case of (NUM) problem, the update of
the link price $\hat{\lambda}_j^k$ requires only the source rates
$z_i^k$ which use link $j$. Note that analysis of step $3$ in the
Algorithm (\textbf{DFG}) can be derived in a similar way as for
step $2$. Also, step $2$ in phases 1 and 2 and step 3 in phase 1
of the Algorithm (\textbf{H-DFG}) follows similarly. Note also
that Algorithm (\textbf{DG}) has the same iterations as phase 2 of
Algorithm (\textbf{H-DFG}).

Further, we note that all the estimates for the convergence rate
for primal and dual suboptimality and primal feasibility violation
derived in Sections \ref{subsec_convergence_dfwg},
\ref{subsec_convergence_h-dfg} and
\ref{subsec_convergence_dfg_error_bound} depends on the upper
bound on the norm of the optimal Lagrange multipliers
$\mathcal{R}$, which at its turn depends on the degree of
separability of problem \eqref{eq_prob_princc}, characterized by
the sets $\Ncal_i$ and $\bar{\Ncal}_j$. In order to see this
dependence we can write further:
\begin{equation}
\label{eq_dependce_sparsity} \mathcal{R}^2=\max_{\lambda^* \in
\Lambda^*} \|\lambda^*\|_W^2=\max_{\lambda^* \in \Lambda^*}
\sum_{j=1}^{\bar{M}} \sum_{i \in \bar{\Ncal}_j}L_{\mathrm{d}_i}
\|\lambda^*_j\|^2,
\end{equation}
from which it is straightforward to notice that $\mathcal{R}$
depends on the cardinality of each $\bar{\Ncal}_j$. On the other
hand, for each $i$ we recall that:
\begin{equation*}
L_{\text{d}_i}=\frac{\left\|\left[
\begin{array}{c}
  \left[A_{ji}\right]_{j \in
\Ncal_i} \\
   \left[C_{ji}\right]_{j \in
\Ncal_i} \\
\end{array}
\right]\right\|^2}{\sigma_i},
\end{equation*}
which depends on the cardinality of the set $\Ncal_i$. Thus, we
can conclude that $\mathcal{R}$ depends on the cardinality of
$\Ncal_i$ and $\bar{\Ncal}_j$ which represent a natural measure
for the degree of separability of our original problem
\eqref{eq_prob_princc}.

\section{Numerical simulations}
\label{sec_numerical}

In order to certify the theoretical results previously presented, in this section we test the performances of Algorithms
(\textbf{DFG}), (\textbf{H-DFG}) and distributed dual gradient Algorithm (\textbf{DG}) for
solving the (DC-OPF) problem in form \eqref{eq_prob_opf} for different IEEE bus test cases. We recall that in the Algorithm (\textbf{DG}),
at each iteration $k$ the dual variable is updated as follows:
\begin{equation*}
\lambda^{k+1}=\left[\lambda^k+W^{-1} \nabla d(\lambda^k)\right]_{\mathbb{D}}.
\end{equation*}
The numerical simulation are performed on different power systems, representing classical test cases from the
literature \cite{ZimMur:11}, with the number of buses $M$ ranging from 9 to 300, the number of generators from 3 to 69 and the number of
interconnecting lines from 18 to 411.
The descriptions of the power systems
are listed in the table below:
\begin{table}[ht!]
\centering 
\begin{tabular}{|c|c|c|c|c|c|c|c|}
\hline
Test case & $M$ & $M_g$  & $\bar{M}$ & $n$ & $p$ & $q$ & Details\\
\hline
 pws$_1$  &  9  &  3  &  9  & 12  & 9  &  18  &  Example taken from \cite{ZimMur:11}  \\
\hline
 pws$_2$   &  14  &  5  &  20  & 19  & 14 &  40  &  IEEE 14 bus test case  \\
\hline
 pws$_3$   &  30  &  6  &  41  & 36  & 30 &  82  &  IEEE 30 bus test case  \\
\hline
 pws$_4$   &  39  &  10  &  46  & 49  & 39  &  92  & 39 bus New England system   \\
\hline
 pws$_5$   &  57  &  7  &  80  & 64  & 57  &  160  &  IEEE 57 bus test case  \\
\hline
 pws$_6$   &  118  &  54  &  186  & 172  & 118  &  372  &  IEEE 118 bus test case  \\
\hline
 pws$_7$   &  300  &  69  &  411  & 369  & 300  &  822  &  IEEE 300 bus test case  \\
 \hline
\end{tabular}
\caption{Description of the test cases.}
\label{table_test_cases}
\end{table}

For each power system considered for simulation we generate the local constraints sets imposed on the phase angle and on the generated power of each bus $i$, $\Theta_i$ and $\mathcal{P}_i$, respectively, the local loads $P_i^d$ and the matrices $E$, $R$ and $A^g$ using the data extract from the MATPOWER toolbox \cite{ZimMur:11}. Also, for each test case we take the parameters of the local cost functions as follows: $q_i=2$, $p_i=10$, $\gamma_i=2$ and $\beta_i=0.1$.

In the case of (DC-OPF) problem the Lagrangian function takes the following form:
\begin{align*}
\lu(\theta,P^g,\lambda)=\sum_{i=1}^M f_i(\theta_i,P_i^g)+\langle \nu, E^TRE \theta -A^g P^g +P^d \rangle
+\langle \mu, \left[\begin{array}{c}
                       ~~RE \\
                       -RE
                     \end{array}
\right] \theta -\left[\begin{array}{c}
                       -\overline{F} \\
                       \underline{F}
                     \end{array}
\right] \rangle
\end{align*}
where we recall that $\theta=\left[\theta_1^T \cdots \theta_M^T\right]^T$, $P^g=\left[\left(P_1^g\right)^T \cdots \left(P_{\bar{M}}^g\right)^T\right]^T$, the functions $f_i$ are given by \eqref{eq_cost_local} if the bus $i$ is directly coupled to a generator unit or by \eqref{eq_cost_local_simple} otherwise and $\lambda=\left[\nu^T~\mu^T\right]^T$. For all algorithms, for each Lagrange multiplier $\lambda$ we have to compute the optimal solution of the inner problem, i.e. the minimization of the Lagrange function subject to the local constraints $\theta_i \in \Theta_i$ and $P_i^g \in \mathcal{P}_i$. As we have shown in Section \ref{subsec_distributed_dfwg}, due to the separability of Lagrangian $\lu$, this can be done distributively, i.e. computing the phase angle $\theta_i(\lambda)$ and the generated power $P^g_i(\lambda)$, for a given trading price $\lambda$, require only local information. Moreover, in the case of (DC-OPF) problem \eqref{eq_prob_opf}, $\theta_i(\lambda)$ and $P^g_i(\lambda)$ can be computed in closed form by solving the following scalar equations derived from the optimality conditions of the inner problems:
\begin{equation*}
\left\{\begin{array}{c}
    q_i\left(\theta_i-\theta_i^{\text{ref}}\right)+\sum\limits_{j \in \mathcal{S}_i} \nu_j \left[E^TRE\right]_{ji}+\sum \limits_{l \in \mathcal{N}_i}  \mu_l^T \left[\begin{array}{c}                                                                                                                                                     ~~[RE]_{li}\\
    -[RE]_{li}\end{array}\right]=0
    \\
    p_i\left(P^g_i-P_i^{g,\text{ref}}\right)+ \nu_i A^g_{ij_{i}}-\frac{\gamma_i}{\beta_i+P_i^g}=0,
  \end{array}\right.
\end{equation*}
for all buses $i \in V_1$, where $j_i$ denotes the position of the generator unit in $P^g$ directly coupled to bus $i$. In order to compute $\theta^k_i=\theta_i(\lambda^k)$ and $P^{g^k}_i=P^g_i(\lambda^k)$ for an iteration $k$, after solving the previous equations we have to project their solutions onto the local box constraints sets $\Theta_i$ and $\mathcal{P}_i$.

The reader should note that  in the context of (DC-OPF) problem,
$\nu$ multipliers associated to the power balance equation have the
economic interpretation as the optimal energy trading prices at the buses of the network. Therefore, our algorithms are able to
identify also the  optimal energy pricing rates for the energy
traded through the interconnections in a distributed fashion.
Thus, it is not necessary to set up a common control center, but
it is sufficient to interchange a small amount of information among
the involved buses. Moreover, the update of the trading prices (dual variables) can be also done in a distributed fashion as follows:
\begin{equation*}
\left\{\begin{array}{c}
    \hat{\mu}_l^{k+1}=\mu_l^{k}+W_{\mu_{ll}}^{-1} \sum \limits_{i \in \bar{\Ncal}_l} \left(\left[\begin{array}{c}
                                                                                             ~~\left[RE\right]_{li} \\
                                                                                             -\left[RE\right]_{li}
                                                                                           \end{array}\right]\theta_i^k-\left[\begin{array}{c}
                              \overline{F}_l \\
                              -\underline{F}_l
                            \end{array}
    \right]\right)
    \\
    \hat{\nu}_j^{k+1}=\nu_j^{k}+W_{\nu_{jj}}^{-1} \sum \limits_{i \in \Scal_j} \left(\left[E^TRE\right]_{ji} \theta_i^k-A^g_{ij_{i}}P_i^{g^k}+P^d_j\right),
  \end{array}\right.
\end{equation*}
for all lines $l \in V_2$ and buses $j \in V_1$. We solve the (DC-OPF) problem \eqref{eq_prob_opf} using Algorithms (\textbf{DFG}), (\textbf{H-DFG}) and (\textbf{DG}) and we compare their performances in terms of the number of iterations. We also consider the centralized versions of these algorithms, namely: (\textbf{CFG}), (\textbf{H-CFG}) and (\textbf{CG}), where by centralized version we understand the version of the algorithm where instead of the step size given by  matrix $W$ we use $L_{\text{d}}I_{p+q}$ with $L_{\text{d}}$ denoting the Lipschitz constant of the gradient $\nabla d$ of the dual function. We recall that the optimization variable $z=\left[z_1^T \cdots z_M^T\right]^T$, where $z_i=\theta_i$ for the buses which do not have a generator unit and $z_i=\left[\theta_i^T ~P_i^{gT}\right]^T$ for the ones directly coupled to a generator. In order to construct the matrix $A$ we interpolate the columns of $E^TRE$ and $A^g$ on the corresponding positions, while $C$ is formed by intercalating in the matrix $\left[(RE)^T~ -(RE)^T\right]^T$ columns with elements equal to zero on the positions corresponding to the position of $P_i^g$ in the vector $x$.

In Table \ref{table_iterations_test_cases} we show, for each test case, the number of iterations performed by the algorithms in order to find a suboptimal primal solution $\hat{z}^k$ which satisfy the following stopping criteria for primal suboptimality and feasibility violation:
\begin{equation}
\label{eq_criteria_table}
\frac{|f(\hat{z}^k)-f^*|}{f^*} \leq \epsilon~\text{and}~\left\|\left[G \hat{z}^k-g
\right]_{\mathbb{D}}\right\|_{W^{-1}} \leq \epsilon,
\end{equation}
where we recall that $G=\left[A^T ~C^T\right]^T$ and $g=\left[b^T ~c^T\right]^T$. Note that in the case of Algorithm (\textbf{DFG}) $\hat{z}^k$ is given by \eqref{primal_point}, while for Algorithms (\textbf{H-DFG}) and (\textbf{DG}) $\hat{z}^k=z^{k^*}$ and $\hat{z}^k=z(\lambda^k)$, respectively. We also consider the same estimates for the centralized version of the algorithms. In our simulation we consider an accuracy $\epsilon=0.01$. It is straightforward to notice that the suboptimality criterion satisfied with this accuracy implies the fact that the difference between the value of the cost function $f(\hat{z}^k)$ and the optimal value $f^*$ is less than $1\%$. For each test case, we use CVX in order to compute the optimal value $f^*$. Also, in the case when the imposed accuracy has not been attained after $3 \cdot 10^5$ iterations, we stoped the algorithm and reported $*$.

\begin{table*}[ht!]
\centering {  
\begin{tabular}{|c|c|c|c|c|c|c|}
\hline
\backslashbox{Test case}{Algorithm} & \textbf{DFG} & \textbf{CFG} & \textbf{H-DFG} & \textbf{H-CFG} & \textbf{DG} & \textbf{CG} \\
\hline
 pws$_1$   &  4486  &  4134  &  700  &  646  &  168619  &  143283  \\
\hline
 pws$_2$  &  1991  &  1920  &  944  &  1066  &  203210  &  214746  \\
\hline
 pws$_3$  &  1368  & 2013   & 503   &  1356  &  27026  &  52893  \\
\hline
 pws$_4$   &  1756  &  6343  &  1316  &  4835  &  69961  &  275343  \\
\hline
 pws$_5$   &  4876  & 21123  &  2003  &  15507  &  $*$  &  $*$  \\
\hline
 pws$_6$   &  8117  &  45787  &  5787  &  35624  &  $*$  &  $*$  \\
\hline
pws$_7$   &  19432  &  63456  &  9978  &  67843  &  $*$  &  $*$  \\
\hline
\end{tabular}}
\caption{Number of iterations performed for finding an $\epsilon$-suboptimal solution of the (DC-OPF) problem for each test case.} \label{table_iterations_test_cases}
\end{table*}

Some remarks are worth to be mentioned. First, we can observe from Table \ref{table_iterations_test_cases} that both the proposed Algorithms (\textbf{DFG}) and (\textbf{H-DFG}) clearly outperform the classical dual gradient Algorithm (\textbf{DG}). Thus, the practical behaviour observed in simulations certifies the
theoretical results derived in the previous sections, where we have proved that the rate of convergence of the proposed algorithms improves the well known rate of convergence of order $\mathcal{O}(\frac{1}{k})$ for the Algorithm (\textbf{DG}). This behaviour is also valid for the centralized case. Another important aspect consists in the fact that for all algorithms, when the dimension of the problem increases, the distributed version becomes more efficient than the centralized one. This is a consequence of the fact that when the number of busses increases, the level of sparsity of the matrices $A$ and $C$, characterized in terms of the indices sets $\Scal_i$, $\bar{\Scal}_i$ and $\Ncal_l$, is high and therefore the Lipschitz constants $L_{\text{d}_i}$ are small in comparison with the overall Lipschitz constant $L_{\text{d}}$ (see Section \ref{subsec_distributed_dfwg} for a more detailed discussion). These differences between $L_{\text{d}_i}$ and $L_{\text{d}}$ lead to a grater step size in the case of distributed algorithms in comparison with the centralized ones,  thus the distributed algorithms perform faster. 

\begin{figure}[h!]
\centering
\includegraphics[angle=0,height=5.5cm,width=8.5cm]{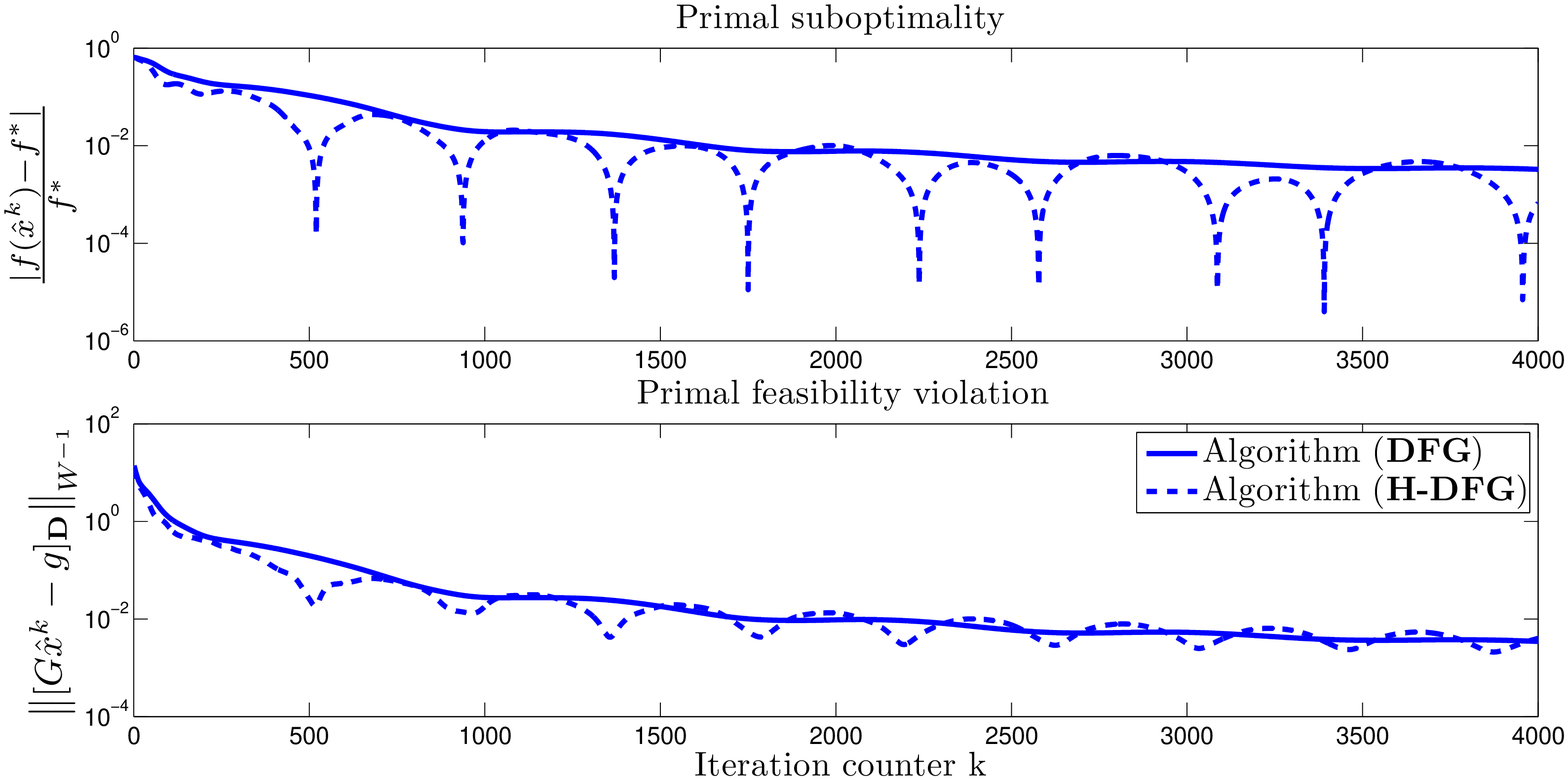}
\caption{Comparison between Algorithms
\textbf{(DFG)} and (\textbf{H-DFG}).}
\label{fig:plot_pws4}
\end{figure}

Further, we are also interested in analyzing the behaviour of the proposed algorithms in terms of the primal suboptimality and feasibility violation. For this purpose we consider the 39 bus New England system (pws$_4$). For this test case we have a number of $M=39$ buses, $M_g=10$ generator units and $\bar{M}=46$ lines between buses. We let the Algorithms (\textbf{DFG}) and (\textbf{H-DFG}) perform a number of 4000 iterations and we show in Figure \ref{fig:plot_pws4} the evolution of primal suboptimality and feasibility violation for each algorithm.
\vspace*{-0.3cm}

We can observe that, on the one hand, the Algorithm (\textbf{H-DFG}) is faster than (\textbf{DFG}) but on the other hand both primal suboptimality and primal feasibility for Algorithm (\textbf{H-DFG}) have an oscillating behaviour, while in the case of Algorithm (\textbf{DFG}) these quantities have a smooth evolution.

\begin{figure}[h!]
\centering
\includegraphics[angle=0,height=5.5cm,width=8.5cm]{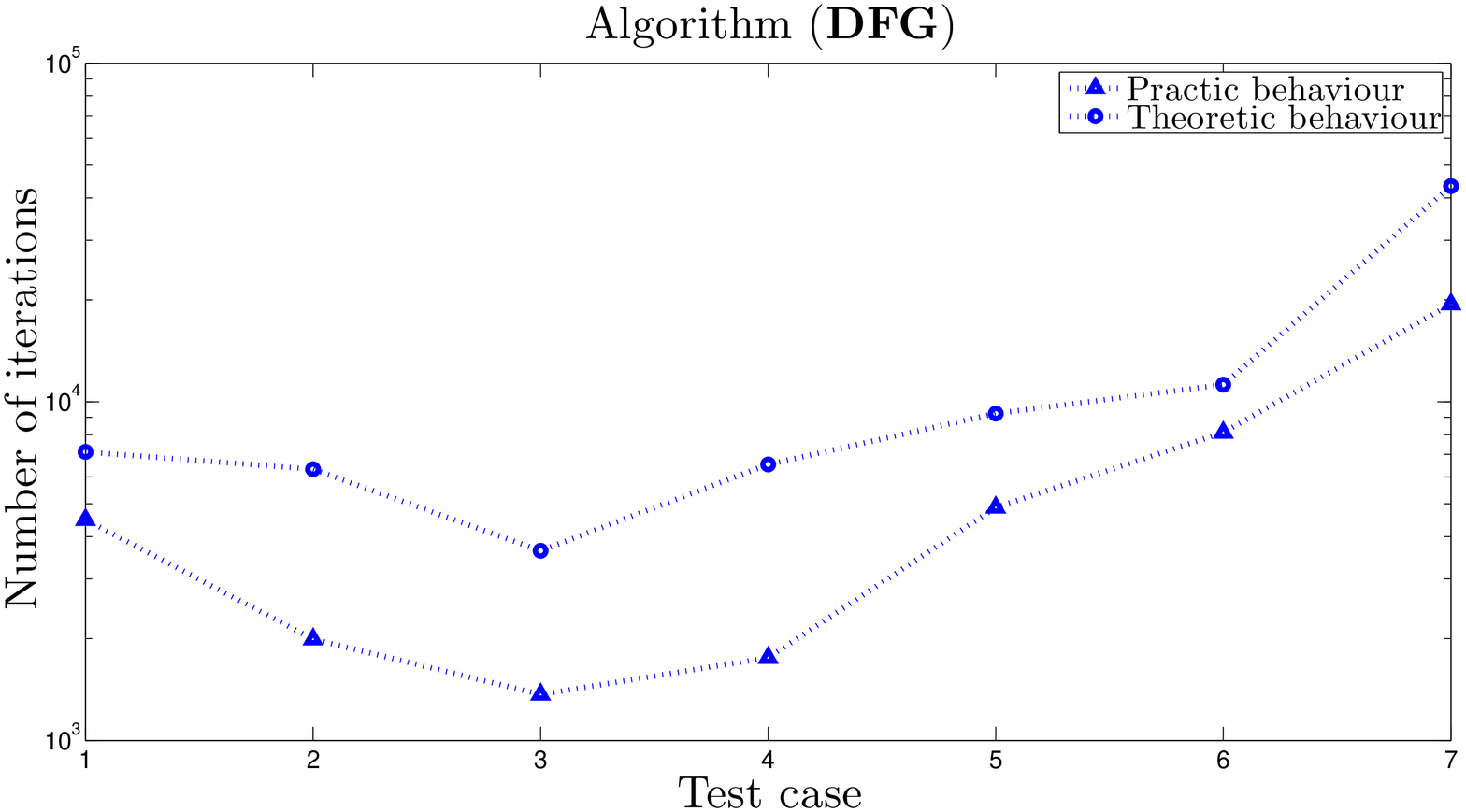}
\caption{Comparison between the real number of iterations performed by Algorithm
\textbf{(DFG)} in practice and the theoretic number of iterations.}
\label{fig:plot_teoretic}
\end{figure}

For Algorithm (\textbf{DFG}) we also plot in Figure \ref{fig:plot_teoretic} the real number of iterations observed in practice and the theoretic number
of iterations derived in Section \ref{sec_ddfg}. We can observe from Figure \ref{fig:plot_teoretic} that the estimates obtained for the number of iterations are closed to the real number of iterations performed by the algorithm in practice.


%
%
%
%
%
%
%
%



%


\bibliographystyle{plain}
\bibliography{bibliografie2013}

\begin{thebibliography}{10}

\bibitem{BecNed:13}
A.~Beck, A.~Nedic, A.~Ozdaglar, and M.~Teboulle.
\newblock Optimal distributed gradient methods for network resource allocation
  problems.
\newblock Technical report, Technion, 2013.

\bibitem{DoaKev:11}
M.D. Doan, T.~Keviczky, and B.~De Schutter.
\newblock A distributed optimization-based approach for hierarchical mpc of
  large-scale systems with coupled dynamics and constraints.
\newblock In {\em Proceedings of 50th IEEE Conference on Decision and Control},
  pages 5236--5241, 2011.

\bibitem{Gau:77}
J.~Gauvin.
\newblock A necessary and sufficient regularity condition to have bounded
  multipliers in nonconvex programming.
\newblock {\em Mathematical Programming}, 12:136--138, 1977.

\bibitem{KiwLar:07}
K.C. Kiwiel, T.~Larsson, and P.O. Lindberg.
\newblock Lagrangian relaxation via ballstep subgradient methods.
\newblock {\em Mathematics of Operations Research}, 32(3):669--686, 2007.

\bibitem{KogFin:11}
M.~Kogel and R.~Findeisen.
\newblock Fast predictive control of linear systems combining nesterov's
  gradient method and the method of multipliers.
\newblock In {\em Proceedings of 50th IEEE Conference on Decision and Control},
  pages 501--506, 2011.

\bibitem{LarPat:98}
T.~Larsson, M.~Patriksson, and A.~Stromberg.
\newblock Ergodic convergence in subgradient optimization.
\newblock {\em Optimization Methods and Software}, 9(1--3):93--120, 1998.

\bibitem{LuoTse:92a}
Z.Q. Luo and P.~Tseng.
\newblock On the convergence of coordinate descent method for convex
  differentiable minimization.
\newblock {\em Journal of Optimization Theory and Applications}, 72(1):7--35,
  1992.

\bibitem{LuoTse:93a}
Z.Q. Luo and P.~Tseng.
\newblock On the convergence rate of dual ascent methods for linearly
  constrained convex minimization.
\newblock {\em Mathematics of Operations Research}, 18(4):846--867, 1993.

\bibitem{Man:85}
O.L. Mangasarian.
\newblock Computable numerical bounds for lagrange multipliers of stationary
  points of non-convex differentiable non-linear programs.
\newblock {\em Operations Research Letters}, 4(2):47--48, 1985.

\bibitem{NecNed:13}
I.~Necoara and V.~Nedelcu.
\newblock Rate analysis of inexact dual first order methods: application to
  dual decomposition.
\newblock {\em IEEE Transactions on Automatic Control}, accepted, 2013.

\bibitem{NecNed:14}
I.~Necoara and V.~Nedelcu.
\newblock Distributed dc optimal power flow based on dual fast gradient
  methods.
\newblock {\em IEEE Transactions on Power Systems}, submitted, 2014.

\bibitem{NecNes:11}
I.~Necoara, Y.~Nesterov, and F.~Glineur.
\newblock A random coordinate descent method on large optimization problems
  with linear constraints.
\newblock Technical report, University Politehnica Bucharest, 2011.

\bibitem{NecSuy:08}
I.~Necoara and J.A.K. Suykens.
\newblock Application of a smoothing technique to decomposition in convex
  optimization.
\newblock {\em IEEE Transactions on Automatic Control}, 53(11):2674--2679,
  2008.

\bibitem{Ned:13}
V.~Nedelcu.
\newblock Rate analysis of dual gradient methods. application to control
  problems.
\newblock Technical report, University Politehnica of Bucharest, Bucharest,
  Romania, december 2013.
\newblock {\verb http://141.85.225.150/papers/nedelcu.pdf }.

\bibitem{NedOzd:09}
A.~Nedic and A.~Ozdaglar.
\newblock Approximate primal solutions and rate analysis for dual subgradient
  methods.
\newblock {\em SIAM Journal on Optimization}, 19(4):1757--1780, 2009.

\bibitem{Nes:04}
Y.~Nesterov.
\newblock {\em Introductory Lectures on Convex Optimization: A Basic Course}.
\newblock Kluwer, Boston, USA, 2004.

\bibitem{Nes:05}
Y.~Nesterov.
\newblock Smooth minimization of non-smooth functions.
\newblock {\em Mathematical Programming}, 103(1):127--152, 2005.

\bibitem{Nes:12a}
Y.~Nesterov.
\newblock How to make the gradients small.
\newblock {\em Optima}, 88:10--11, 2012.

\bibitem{PatBem:12}
P.~Patrinos and A.~Bemporad.
\newblock An accelerated dual gradient-projection algorithm for embedded linear
  model predictive control.
\newblock {\em IEEE Transactions on Automatic Control}, under review, 2012.

\bibitem{RicMor:11}
S.~Richter, M.~Morari, and C.N. Jones.
\newblock Towards computational complexity certification for constrained mpc
  based on lagrange relaxation and the fast gradient method.
\newblock In {\em Proceedings of 50th IEEE Conference on Decision and Control},
  pages 5223--5229, 2011.

\bibitem{Rob:73}
S.~M. Robinson.
\newblock Bounds for error in the solution set of a perturbed linear program.
\newblock {\em Linear Algebra and its Applications}, 6:69--81, 1973.

\bibitem{RocWet:98}
R.T. Rockafellar and R.J. Wets.
\newblock {\em Variational Analysis}.
\newblock Springer-Verlag, New York, 1998.

\bibitem{ScoMay:99}
P.O.M. Scokaert, D.Q. Mayne, and J.B. Rawlings.
\newblock Suboptimal model predictive control (feasibility implies stability).
\newblock {\em IEEE Transactions on Automatic Control}, 44(3):648--654, 1999.

\bibitem{SenShe:86}
S.~Sen and H.D. Sherali.
\newblock A class of convergent primal-dual subgradient algorithms for
  decomposable convex programs.
\newblock {\em Mathematical Programming}, 35(3):279--297, 1986.

\bibitem{TriZym:08}
N.~Trichakis, A.~Zymnis, and S.~Boyd.
\newblock Utility maximization with delivery contracts.
\newblock In {\em Proceedings of 17th IFAC World Congress}, 2008.

\bibitem{WanLin:13}
P.W. Wang and C.J. Lin.
\newblock Iteration complexity of feasible descent methods for convex
  optimization.
\newblock Technical report, Department of Computer Science, National Taiwan
  University, 2013.

\bibitem{XiaBoy:06}
L.~Xiao and S.~Boyd.
\newblock Optimal scaling of a gradient method for distributed resource
  allocation.
\newblock {\em Journal of Optimization Theory and Applications}, 129(3), 2006.

\bibitem{ZimMur:11}
R.D. Zimmerman, C.E. Murillo-Sanchez, and R.J. Thomas.
\newblock Matpower: Steady-state operations, planning, and analysis tools for
  power systems research and education.
\newblock {\em IEEE Transactions on Power Systems}, 26(1):12--19, 2011.

\end{thebibliography}


\end{document}